\documentclass[11pt]{amsart}
\usepackage{eucal}
\usepackage{tikz}
\usepackage{amsfonts,amssymb}
\usepackage{epsfig}
\usepackage{enumerate}

\newtheorem{theorem}{Theorem}[section]
\newtheorem{lemma}[theorem]{Lemma}

\newtheorem{corollary}[theorem]{Corollary}
\newtheorem{defn}[theorem]{Definition}

\newcommand{\Q}{\mathbb Q}
\newcommand{\C}{\mathbb C}
\newcommand{\R}{\mathbb R}

\renewcommand{\P}{\mathbb P}

\newcommand{\Z}{\mathbb Z}
\newcommand{\N}{\mathbb N}

\newcommand{\bracket}[1]{\{#1\}}

\title{On the symplectic invariance of log Kodaira dimension}
\author{Mark McLean}

\begin{document}

\begin{abstract}
Suppose that $A$ and $B$ are symplectomorphic smooth affine varieties.
If $A$ is acylic of dimension $2$ then $B$ has the same log Kodaira dimension as $A$.
If the dimension of $A$ is $3$, has log Kodaira dimension $2$ and satisfies some other conditions
then $B$ cannot be of log general type.
We also show that if $A$ and $B$ are symplectomorphic affine varieties of any dimension
then any compactification of $A$ by a projective variety is uniruled
if and only if any such compactification of $B$ is uniruled.
\end{abstract}

\maketitle

\bibliographystyle{halpha}


\tableofcontents

\section{Introduction}
Much of algebraic geometry is governed by the numerical properties of the canonical class.
Other useful properties such as uniruledness and rational connectivity
have also played a major role.
Also closed smooth projective varieties, and more generally K\"{a}hler manifolds
have been studied extensively from both a topological and symplectic viewpoint.
In complex dimension $2$, tools such as Donaldson theory and Seiberg Witten invariants have
been used to study such manifolds.
For instance in \cite{Witten:Monopoles} algebraic surfaces of general type
have plus or minus their canonical class as diffeomorphism invariants.
Plurigenera and hence Kodaira dimension for algebraic surfaces
are also shown in \cite{FriedmanMorgan:AlgebraicSW}
to be diffeomorphism invariants using Seiberg Witten theory.
There are many other results for these surfaces.

One can also study these varieties from a symplectic perspective.
We can use tools such as Gromov-Witten
theory to see what the symplectic structure tells us about our algebraic variety.
%
Work by Koll\'{a}r and Ruan (\cite{Kollar:rationalcurves} and \cite{Ruan:virtual})
tells us the property of being uniruled is a symplectic invariant using Gromov Witten theory.
Another extremely useful notion in algebraic geometry is rational connectedness
and this has been studied from a symplectic viewpoint in  \cite{Voisin:rationallyconnected} and
\cite{ZTian:rationallyconnected3fold}.

Less has been done to study open algebraic varieties
from a symplectic perspective, although there has been some work
\cite{TianJunZhang:additivity}.
Also there isn't as much work in higher dimensions
although there is some progress in dimension $3$
(see for instance \cite{Ruan:3folds}).
This paper addresses some of these issues. We will be primarily concerned with smooth affine
varieties and we will study them from a symplectic perspective.
Every smooth affine variety has a symplectic structure coming from some embedding in $\C^N$
and this is a biholomorphic invariant (see \cite{EliahbergGromov:convexsymplecticmanifolds}).
A particular algebraic invariant of $A$
is called the log Kodaira dimension. One can ask
to what extent is the log Kodaira dimension a symplectic invariant?
Log Kodaira dimension is a number $\overline{\kappa}(A)$ which takes values in
$\{-\infty,0,1,\cdots,\text{dim}_\C A\}$.
We say that $A$ is of log general type if $\overline{\kappa}(A) = \text{dim}_\C A$.
A precise definition is given at the start of section \ref{section:logkodaira}.

We show that log Kodaira dimension is a symplectic invariant for smooth acyclic
affine surfaces (Theorem \ref{theorem:acyclicsurfaceinvariance}).
We also show that if $A$ and $B$ are symplectomorphic smooth affine varieties  such that:
\begin{enumerate}
\item $A$ has complex dimension $3$.
\item $A$ can be compactified by a smooth normal crossing nef divisor
which is linearly equivalent to some smooth divisor
\item The log Kodaira dimension of $A$ is $2$.
\end{enumerate}
then the log Kodaira dimension of $B$ is $\leq 2$
(see Theorem \ref{theorem:logkodairaresultsindimension3}).

A projective variety is {\it uniruled} if there is a rational curve passing through every point.
Let $P$ and $Q$ be smooth projective varieties compactifying smooth affine varieties $A$ and $B$ respectively.
We show that if $A$ is symplectomorphic to $B$ and $P$ is uniruled then $Q$ is also uniruled
(Theorem \ref{theorem:birationalinvarianceofuniruledness}).

In order to prove these theorems we introduce three notions of uniruledness for smooth affine
varieties. The first notion is defined for an object called the Liouville domain associated to our affine variety.
This is a symplectic invariant and is defined in Section \ref{section:uniruledliouville}.
The second notion says that a smooth affine variety is algebraically $k$ uniruled if there is a morphism
from $\P^1$ minus at most  $k$ points to our variety passing through a generic point
and is defined in Section \ref{section:affineuniruled}.
The third notion defined in Section \ref{section:affineuniruled2} is more flexible than the second notion as it now involves $J$ holomorphic
curves from $\P^1$ minus some points where $J$ is any appropriate almost complex structure.
This notion is called compactified $k$ uniruled.
We show by using degeneration to the normal cone techniques that the first definition implies the second definition
(Theorem \ref{theorem:kuniruledimpliesalgebraicallyuniruled}).
Also using other simpler techniques we can show that  third definition implies the first one
(Theorem \ref{theorem:algebraicuniruledimplesuniruled}).
Putting all of this together one gets that if $A$ is symplectomorphic to $B$ and is compactified $k$
uniruled then $B$ is algebraically $k$ uniruled.

Because there is a relationship between log Kodaira dimension and uniruledness
in low dimensions (see \cite{MiyanishiSugie:affinessurfaces},\cite{Miyanishi:openaffine},
\cite{Kawamata:classification} and \cite{Kishimoto:affinethreefolds})
we obtain our log Kodaira dimension results.
Similarly if $P$ and $Q$ are projective with symplectomorphic affine open subsets $A$ and $B$
such that $P$ is uniruled,
then one can show that $A$ is compactified $k$ uniruled for some $k$. Hence $B$ is algebraically $k$
uniruled which in turn implies that $Q$ is uniruled.

The paper is organized as follows:
In Section \ref{section:uniruledliouville} we introduce the reader to uniruled Liouville domains
(first definition).
These are purely symplectic objects.
In Section \ref{section:affineuniruled} we give a purely algebraic definition of uniruledness
for smooth affine varieties (second definition) and relate it to uniruled Liouville domains.
In Section \ref{section:GWintro} we give an introduction to Gromov Witten invariants,
then in Section \ref{section:affineuniruled2} we give a much more flexible
definition of uniruledness (third definition) for smooth affine varieties.
In Section \ref{section:logkodaira} we use all of the above machinery to
prove our log Kodaira dimension invariance results and finally in
Section \ref{section:unirulednesscompactifications} we prove that
projective varieties with symplectomorphic open affine subsets
are either both uniruled or both not uniruled.

\bigskip

{\it Acknowledgements:} I would like to thank Paul Seidel and Ivan Smith for their help in this project.
The author was partially supported by
NSF grant DMS-1005365 and also a grant to the Institute for Advanced Study by
The Fund for Math.

\section{Uniruled Liouville domains} \label{section:uniruledliouville}

Throughout this paper we will use the following notation.
If $U$ is any subset of a topological space then we write $U^o$ for the interior of $U$.
Also if $(N,\omega)$ is a symplectic manifold and $\theta$ is a $1$-form
then we write $X_{\theta}$ to be the unique vector satisfying $\iota_{X_{\theta}} \omega = \theta$.

Let $M$ be a compact manifold with boundary with a $1$-form $\theta_M$ satisfying:
\begin{enumerate}
\item $\omega_M := d\theta_M$ is a symplectic form.
\item The $\omega_M$-dual $X_{\theta_M}$ of $\theta_M$ points outwards along $\partial M$.
\end{enumerate}
We say that $(M,\theta_M)$ is a {\it Liouville domain} if it satisfies the above properties.
Let $J$ be an almost complex structure compatible with the symplectic form $\omega_M$.
We say that $J$ is a {\it convex almost complex structure} on $M$
if there is some function $\phi : M \rightarrow \R$ so that:
\begin{enumerate}
\item $\partial M$ a regular level set of $\phi$ and
$\phi$ attains its maximum on $\partial M$.
\item $\theta_M \circ J = d\phi$ near $\partial M$.
\end{enumerate}
Suppose that $(N, \omega_N)$ is a symplectic manifold and let $J_N$ be an almost complex structure.
If $u : S \rightarrow N$ is a $J_N$-holomorphic map from a Riemann surface $S$ to $N$
then the {\it energy of} $u$ is defined to be $\int_S u^* \omega_N$.

\begin{defn}
Let $k>0$ be an integer and $\lambda > 0$ a real number. We say that a Liouville domain $M$
is {\bf $(k,\Lambda)$-uniruled} if for every convex almost complex structure $J$
on $M$ and every point
$p \in M^o$ where $J$ is integrable on a neighbourhood of $p$,
there is a proper $J$ holomorphic map $u : S \rightarrow M^o$ to the interior $M^o$ of $M$
passing through this point.
We require that $S$ is a genus zero Riemann surface, the rank of $H_1(S,\Q)$ is at most $k-1$
and the energy of $u$ is at most $\Lambda$.
\end{defn}

\begin{theorem} \label{theorem:subdomainuniruled}
Suppose that $N,M$ are Liouville domains such that $M$ is a codimension $0$ symplectic submanifold of $N$
with the property that there exists some $1$-form $\theta'$ on $N$ so that $\theta'|_M - \theta_M$ is exact
and so that $d\theta' = d\theta_N$.
If $N$ is $(k,\Lambda)$-uniruled then $M$ is also $(k,\Lambda)$-uniruled.
In particular, the above fact is true if $M$ is a codimension $0$ exact submanifold of $N$
or if the inclusion map $M \hookrightarrow N$ is a symplectic embedding and a homotopy equivalence.
\end{theorem}
Before we prove this theorem we need some preliminary lemmas and definitions.
The following definitions are technically not relevant for the theorem above,
but one of the lemmas used in proving this theorem will also be used later
on in a slightly more general context.
A {\it nodal Riemann surface} is a $1$ dimensional complex analytic variety
with the property that the only singularities are nodal.
We say that it has {\it arithmetic genus $0$} if it can be holomorphically embedded
into a simply connected compact nodal Riemann surface.
An example of an arithmetic genus zero surface is:
$B(1) \cap \{z_1z_2 = 0\} \subset \C^2$ where $B(1)$ is the open unit ball
and $z_1,z_2$ are coordinates for $\C^2$.
Note that a genus zero nodal Riemann surface is a union $S_1,\cdots,S_k$
of smooth Riemann surfaces which only intersect each other at the nodal singularities of $S$.
Here $S_1,\cdots,S_k$ are called the {\it irreducible components of $S$}.
An {\it arithmetic genus $0$ nodal Riemann surface with boundary}
is a closed subset $S$
of a compact arithmetic genus $0$ nodal Riemann surface $C$ with the property that
away from the nodes of $C$, $S$ is a Riemann surface with boundary.
We require that the closure of this boundary does not intersect the nodes of $C$.
This means that the boundary is a union of circles.
An example of such a holomorphic object would be
the closure of $B(1) \cap \{z_1z_2 = 0\}$.
Again an arithmetic genus $0$ nodal Riemann surface with boundary is a union
of smooth Riemann surfaces with boundary intersecting each other
at the nodal singularities away from their boundaries.
These smooth Riemann surfaces with boundary are called the {\it irreducible components of $S$}.
We can form a graph $\Gamma_S$ whose nodes are the irreducible components of $S$
and if two irreducible components intersect at some point $p$
then we have an edge $E_p$ joining the appropriate nodes.
This is called the {\it dual graph} of $S$.
The dual graph of every connected arithmetic genus zero compact Riemann surface
is a tree.

\begin{tikzpicture} [scale=1.0]
%
\draw[color=lightgray] (5,5) ellipse (1 and 0.25);
\draw (5,5) node {\tiny A} circle (1);
\draw[shift={(5,5)},rotate=35,shift={(-5,-5)}] (7,5)+(0,1) arc (90:270:1);
\draw[shift={(5,5)},rotate=130,shift={(-5,-5)}] (7,5)+(0,1) arc (90:270:1);
\draw[shift={(5,5)},rotate=35,shift={(-5,-5)}] (7,5) node {\tiny B} ellipse (0.25 and 1);
\draw[shift={(5,5)},rotate=130,shift={(-5,-5)}] (7,5) node {\tiny C} ellipse (0.25 and 1);
\draw[shift={(5+5,5)},rotate=130,shift={(-5,-5)}]  (5,5) -- (7,5) node[below,left] {\tiny C};
\draw[shift={(5+5,5)},rotate=35,shift={(-5,-5)}] (5,5) node[below] {\tiny A} -- (7,5) node[below,right] (7,5) {\tiny B};
\draw[fill=black] (5+5,5) circle (0.05);
\draw[shift={(5+5,5)},rotate=35,shift={(-5,-5)},fill=black] (7,5) circle (0.05);
\draw[shift={(5+5,5)},rotate=130,shift={(-5,-5)},fill=black] (7,5) circle (0.05);
\node  at (4,3.5) {  Riemann surface.};
\node  at (4+5,4.5) { Dual graph. };

\end{tikzpicture}

\begin{lemma} \label{lemma:topologyofcurve}
Let $M$ be a Liouville domain and
$N$ any exact symplectic manifold so that $M$ is a codimension $0$
symplectic submanifold of $N$ with the additional property that there is some $1$-form
$\theta'$ on $N$ so that $\theta'|_M - \theta_M$ is exact
and so that $d\theta' = d\theta_N$.
Let $J$ be a compatible almost complex structure on $N$ so that $J$ restricted to $M$
is a convex almost complex structure. If $u : S \rightarrow N^o$ is a $J$-holomorphic curve
with the property that $u^{-1}(M)$ is compact and $S$ is an arithmetic genus zero nodal Riemann surface then
the map $H_1(u^{-1}(M^o)) \rightarrow H_1(S)$ is injective.
\end{lemma}
\proof of Lemma \ref{lemma:topologyofcurve}.
By definition we can choose a collar neighbourhood $(1-\epsilon,1] \times \partial M$ of $\partial M$
inside $M$ so that $d\theta_M \circ J= dr$ where $r$ parameterizes the interval $(1-\epsilon,1]$.
For $R \in (1-\epsilon,1)$ we define $M_R$ to be $M \setminus \{r > R\}$.
We will show that $H_1(u^{-1}(M_R)) \rightarrow H_1(S)$ is injective for generic $R$ and this will prove the
theorem because $H_1(u^{-1}(M^o))$ is the direct limit of $H_1(u^{-1}(M_R))$ as $R$ tends to $1$.

For generic $R$, $\partial M_R$ is transverse to $u$. This means that $S_R := u^{-1}(M_R)$
is an arithmetic genus $0$ compact nodal Riemann surface with boundary.
Also the closure of $S \setminus S_R$ is a possibly non-compact nodal Riemann surface with boundary equal
to $\partial S_R$. We will write this closure as $S^c_R$.
We let $\theta$ be a $1$-form on $N$ so that $\theta'-\theta$ is exact and so that
$\theta = \theta_M$ on a neighbourhood of $M_R$.
The maximum principle \cite[Lemma 7.2]{SeidelAbouzaid:viterbo} using the $1$-form $\theta$ tells us that
every irreducible component of $S^c_R$ is non-compact.
Let $S'_1,\cdots,S'_{l'}$ be these irreducible components.
These are non-compact Riemann surfaces with compact boundary.
Hence they have the property that
$H_1(\partial S'_i) \to H_1(S'_i)$ is injective.
This in turn implies that $H_1(\partial S^c_R) \to H_1(S^c_R)$ is injective.
Because $S$ is the union of $S_R$ and $S^c_R$ along $\partial S^c_R$, we have by
a Mayor-Vietoris argument that the map $H_1(S_R) \rightarrow H_1(S)$ is injective.
Hence $H_1(u^{-1}(M^o)) \rightarrow H_1(S)$ is injective.
\qed

\begin{lemma} \label{lemma:convexalmostcomplexstructureexistence}
Every Liouville domain $(M,\theta_M)$ has a convex almost complex structure.
\end{lemma}
\proof of Lemma \ref{lemma:convexalmostcomplexstructureexistence}.
By flowing $\partial M$ backwards along $X_{\theta_M}$
we have a neighbourhood $(1-\epsilon,1] \times \partial M$ of $M$
so that $\theta_M = r \alpha_M$ where $r$ parameterizes $(1-\epsilon,1]$
and $\alpha_M$ is the contact form $\theta_M|_{\partial M}$.
We define the vector bundle $V_1$ to be the span of the vectors $\frac{\partial}{\partial r}$
and $X_{\theta_M}$ and $V_2$ to be the set of vectors in the kernal of $dr$
and $\theta_M$.
Near $\partial M$, we have that $TM = V_1 \oplus V_2$ and $V_1$, $V_2$ are symplectically orthogonal.
We define $J$ so that:
\begin{enumerate}
\item $J$ is compatible with the symplectic form $\omega_M$.
\item $J(V_1) = V_1$ and $J(V_2) = V_2$ near $\partial M$.
This can be done because these vector spaces are $\omega_M$ orthogonal.
\item $J(r\frac{\partial}{\partial r}) = X_{\theta_M}$.
\end{enumerate}
Here $\theta_M \circ J = dr$ near $\partial M$ and so $J$ is a convex almost complex structure.
\qed

\proof of Theorem \ref{theorem:subdomainuniruled}.
Let $J$ be a convex almost complex structure on $M$.
By Lemma \ref{lemma:convexalmostcomplexstructureexistence}, we have that $N$
admits some convex almost complex structure $J_N$.
Because the space of all almost complex structures compatible with a symplectic form is contractible we
can choose a compatible almost complex structure $J'$ so that $J' = J_N$ near $\partial N$
and $J'|_M = J$.
Let $p \in M$ be a point in the interior of $M$ such that $J$
is integrable on a neighbourhood of $p$.
Because $N$ is $(k,\Lambda)$-uniruled, we have that there exists a proper holomorphic map
$u : S \rightarrow N^o$ passing through $p$ of energy at most $\Lambda$.
Also the rank of $H_1(S,\Q)$ is at most $k-1$.
By Lemma \ref{lemma:topologyofcurve}, the rank of $H_1(u^{-1}(N^o),\Q)$ is also at most $k-1$.
Hence $u|_{u^{-1}(N^o)}$ is a proper $J$ holomorphic map in $M^o$ of energy at most $\Lambda$
and passing through $p$, where $|H_1(u^{-1}(N^o),\Q)| \leq k-1$.
This implies that $M$ is $(k,\Lambda)$-uniruled.
\qed

If $(M,\theta_M)$ is a Liouville domain, then by flowing $\partial M$ backwards along $X_{\theta_M}$
we have a neighbourhood $(1-\epsilon,1] \times \partial M$ of $M$
so that $\theta_M = r \alpha_M$ where $r$ parameterizes $(1-\epsilon,1]$ and where
$\alpha_M$ is a contact form on $\partial M$.
If we glue $[1,\infty) \times \partial M$ to $M$ along $\partial M$ and extend $\theta_M$ by $r\alpha_M$,
then we get a new exact symplectic manifold $\widehat{M}$ called the {\it completion of} $M$.

\begin{theorem} \label{theorem:symplecticinvariance}
Let $M$, $N$ be two Liouville domains such that $\widehat{M}$ is symplectomorphic to $\widehat{N}$.
If $M$ is $(k,\Lambda)$-uniruled then there exists a $\Lambda' > 0$ such that $\widehat{N}$
is $(k,\Lambda')$-uniruled.
\end{theorem}
\proof of Theorem \ref{theorem:symplecticinvariance}.
Let $\phi : \widehat{M} \rightarrow \widehat{N}$ by our symplectomorphism.
By \cite[Lemma 1]{eliashberg:symplectichomology} we can assume that $\phi$ is an exact symplectomorphism
which means that $\phi^* \theta_N = \theta_M + df$ for some function $f$.
Let $\Phi_t : \widehat{M} \rightarrow \widehat{M}$ be the time $t$ flow of the vector field $X_{\theta_M}$.
Because $X_{\theta_M}$ is equal to $r \frac{\partial}{\partial r}$ near infinity where $r$ is the cylindrical
coordinate on $\widehat{M}$, we get that for some $T \geq 0$, $\phi^{-1}(N) \subset \Phi_T(M)$.
Because $\Phi_T^* \theta_M = e^T \theta_M$, we get by a rescaling argument that the Liouville domain
$\phi_T(M)$ is $(k,e^T \Lambda)$-uniruled.
Because $N$ is a codimension $0$ exact symplectic submanifold of $\phi_T(M)$,
we have by Lemma \ref{theorem:subdomainuniruled}
that $N$ is $(k,\Lambda')$-uniruled where $\Lambda' = e^T \Lambda$.
\qed

If we have two Liouville domains $(M,\theta_M)$ and $(N,\theta_N)$
then they are {\it Liouville deformation equivalent} if
there is a diffeomorphism $\phi : M \rightarrow N$ and a smooth family of $1$-forms
$\theta^t_M$ ($t \in [0,1]$) on $M$ with the property that:
\begin{enumerate}
\item $\theta^0_M = \theta_M$,
\item $\theta^1_M = \phi^* \theta_N$, and
\item $(M,\theta^t_M)$ is a Liouville domain for each $t$.
\end{enumerate}

\begin{corollary} \label{corollary:deformationinvariance}
Let $M$, $N$ be two Liouville deformation equivalent Liouville domains.
If $M$ is $(k,\Lambda)$-uniruled then there exists a $\Lambda' > 0$ such that $\widehat{N}$
is $(k,\Lambda')$-uniruled.
\end{corollary}
\proof of Corollary \ref{corollary:deformationinvariance}.
We will first show that $\widehat{M}$ is symplectomorphic to $\widehat{N}$
(which is a standard fact) and then we will use Theorem \ref{theorem:symplecticinvariance}.

Let $\theta^t_M$ be our Liouville deformation on $M$.
By construction, we can complete all the Liouville domains $(M,\theta^t_M)$ giving us
a manifold $\widehat{M}$ and a smooth family of $1$-forms $\theta^t_M$
on $\widehat{M}$ by abuse of notation.
The vector field $X_t$ given by the $d\theta^t_M$-dual of $\frac{d}{dt}(\theta^t_M)$
is integrable.
This is because $dr(X_t)$ is less than or equal to some constant times $r$
where $r$ is the cylindrical coordinate on $\widehat{M}$.
The time $1$ flow of this vector field gives us a symplectomorphism from$(\widehat{M},\theta^0_M)$ to $(\widehat{M},\theta^1_M)$.
So because $\theta_M = \theta^0_M$ and $\phi^* \theta_N = \theta^1_M$ we get that
$(\widehat{M},\theta_M)$ is symplectomorphic to $(\widehat{N},\theta_N)$
Hence by Theorem \ref{theorem:symplecticinvariance}, we have that
$N$ is $(k,\Lambda')$-uniruled for some $\Lambda' \geq 0$.
\qed

\section{Uniruled smooth affine varieties} \label{section:affineuniruled}

We say that an affine variety $A$ is {\it algebraically $k$-uniruled} if through every point $p \in A$
there is a polynomial map $S \rightarrow A$
passing through $p$ where $S$ is equal to $\P^1$ with at most $k$ punctures.
We want to relate this algebraic definition of uniruledness with the one in the last section.
In order to do this we need to associate a Liouville domain with $A$.

\begin{defn} \label{defn:canonicasymplecticformassociatedtoA}
Let $A \subset \C^N$ be a smooth affine variety in $\C^N$. Then we define
the $1$-form $\theta_A$ to be equal to $\sum_{i=1}^N \frac{1}{2}r_i^2 d\vartheta_i$
restricted to $A$ where $(r_i,\vartheta_i)$ are polar coordinates for the
$i$th $\C$ factor. We have that $\omega_A := d\theta_A$ is a biholomorphic invariant
\cite{EliahbergGromov:convexsymplecticmanifolds}.
Also for $R \gg 1$, $(B(R) \cap A, \theta_A)$ is a Liouville domain by
\cite[Lemma 2.1]{McLean:affinegrowth} where $B(R)$ is the closed ball of radius $R$.
We will write $(\overline{A},\theta_A)$ for such a Liouville domain and call it a
{\bf Liouville domain associated to } $A$.
\end{defn}

If $A_1,A_2$ are two isomorphic
smooth affine varieties then any Liouville domain associated to $A_1$
is Liouville deformation equivalent to any Liouville domain associated to $A_2$
by Lemma \ref{lemma:affineLiouvilledomain} in the Appendix.
The problem with the symplectic form $\omega_A$ is that it gives $A$ infinite volume.
But we need to compactify $A$ so in order to deal with this we need another symplectic structure
on $A$ which is compatible with the compactification $X$ of $A$.

\begin{defn} \label{defn:symplecticformonaffinevariety}
Let $A$ be a smooth affine variety and $X$ a smooth projective variety
such that $X \setminus A$ is a smooth normal crossing divisor
(an SNC compactification).
Let $L$ be an ample line bundle on $X$ given by an effective divisor $D$ whose support is $X \setminus A$.
From now on such a line bundle will be called a {\bf line bundle associated to an SNC compactification $X$ of $A$}.
Suppose $|\cdot|$ is some metric on $L$
whose curvature form is a positive $(1,1)$ form. Then if $s$ is some section of $L$ such that $s^{-1}(0) = D$
then we define $\phi_{s,|\cdot|} := -\log{|s|}$ and $\theta_{s,|\cdot|} := -d^c \phi_{s,|\cdot|}$. The two form $d\theta_{s,|\cdot|}$
extends to a symplectic form $\omega_{|\cdot|}$ on $X$ (which is independent of $s$ but does depend on $|\cdot|$).
We will say that $\phi_{s,|\cdot|}$ is a {\bf plurisubharmonic function associated to $L$},
$\theta_{s,|\cdot|}$ a {\bf Liouville form associated $L$} and
$\omega_{|\cdot|}$ a {\bf symplectic form on $X$ associated to $L$}.
\end{defn}

The aim of this section is to prove the following:

\begin{theorem} \label{theorem:kuniruledimpliesalgebraicallyuniruled}
Let $A$ be a smooth affine variety and $\overline{A}$ its associated Liouville domain.
Then if $\overline{A}$ is $(k,\Lambda)$-uniruled then $A$ is algebraically $k$-uniruled.
\end{theorem}
We need some preliminary lemmas before we prove this theorem.

\begin{lemma} \label{lemma:algebraicvarietysymplecticforminvariance}
Let $L$ be a line bundle associated to an SNC compactification $X$ of $A$
and let $|\cdot|_1$, $|\cdot|_2$ be two metrics on $L$ whose curvature forms are positive $(1,1)$-forms.
Then $(A,\omega_{|\cdot|_1})$ is symplectomorphic to $(A,\omega_{|\cdot|_2})$.
\end{lemma}
\proof of Lemma \ref{lemma:algebraicvarietysymplecticforminvariance}.
We have a smooth family of symplectic forms
$\omega_t := (1-t)\omega_{|\cdot|_1} + t\omega_{|\cdot|_2}$ on $X$.
By a Moser argument we have a smooth family of symplectomorphisms
$\phi_t : (X,\omega_{|\cdot|_1}) \to (X,\omega_t)$.
Let $D_t := \phi_t^{-1}(D)$ where $D$ is the associated compactification divisor for $A$.
We have that $(X \setminus D_0,\omega_{|\cdot|_1})$
is symplectomorphic to $(X \setminus D_1, \omega_{|\cdot|_1})$
by  \cite[Lemma 5.15]{McLean:affinegrowth}.
Hence $(A,\omega_{|\cdot|_1})$ is symplectomorphic to $(A,\omega_{|\cdot|_2})$.
\qed

Let $(M,\theta_M)$ be an exact symplectic manifold. We say that it is a
{\it finite type Liouville manifold} if there is an exhausting function $f$
(i.e. it is proper and bounded below) with the property that
$df(X_{\theta_M}) > 0$ outside some compact set. Here
$X_{\theta_M}$ is the $\omega_M := d\theta_M$-dual of $\theta_M$.
We say that a finite type Liouville manifold $M$ is {\it strongly bounded}
if there is some compact set ${\mathcal K} \subset M$ and a constant $T>0$
so that for every point $p$ outside ${\mathcal K}$, the time $T$ flow  of $p$
along $X_{\theta_M}$ does not exist. In other words every point
outside ${\mathcal K}$ flows to infinity within some fixed finite time.

\begin{lemma} \label{lemma:boundedembedding}
Let $(M,\theta_M)$ be a strongly bounded finite type Liouville manifold
and let $f$ be its exhausting function. Choose $C \gg 1$ so that
$df(X_{\theta_M}) >0$ when $f \geq C$. Then there is a $\nu > 0$ 
and an exact symplectic embedding of $(M,\theta_M)$ into the Liouville domain
$(f^{-1}(-\infty,C],\nu \theta_M))$. This embedding is also a homotopy equivalence.
\end{lemma}
\proof of Lemma \ref{lemma:boundedembedding}.
We define $M_C := f^{-1}(-\infty,C]$.
Because $(M,\theta_M)$ is strongly bounded and that $df(X_{\theta_M}) >0$ when $f \geq C$,
we have a constant $T$ so that every flowline starting at a point $p$ with $f(p) \geq C$
flows to infinity in time less than $T$.
Let $\phi_t : M \rightarrow M$ be the time $t$ flow of $-X_{\theta_M}$.
We define an embedding $\iota : M \hookrightarrow M_C$ by $\phi_T$.
The reason why $\phi_T$ sends $M$ into $M_C$ is because we know that points outside $M_C$
flow to infinity in time less than $T$.
We have that $e^T \iota^* \theta_M = \theta_M$.
Hence $\iota$ is an exact symplectic embedding of $(M,\theta_M)$ into
$(M_C,\nu \theta_M)$ where $\nu = e^T$.
\qed

\begin{lemma} \label{lemma:affinevarietycontainedinliouvilledomain}
Let $L$ be a line bundle associated to an SNC compactification $X$ of a smooth affine variety $A$
and let $\theta_{s,\|\cdot\|}$ be a Liouville form associated to $L$.
Then there is some Liouville domain $(M,\theta)$ Liouville deformation equivalent to
$(\overline{A},\theta_A)$
and an exact symplectic embedding $(A,\theta_{s,\|\cdot\|}) \hookrightarrow (M,\theta)$
which is also a homotopy equivalence.
\end{lemma}
\proof of Lemma \ref{lemma:affinevarietycontainedinliouvilledomain}.
In order to prove this Lemma we will show that $(A,\theta_{s,\|\cdot\|})$
is a strongly bounded finite type Liouville manifold and then use Lemma
\ref{lemma:boundedembedding} to finish off the proof.

We write $D := s^{-1}(0)$
and $f_A := \phi_{s,\|\cdot\|} = -\log{\|s\|}$ our plurisubharmonic function associated to $L$.
Locally around some point in $D$
we have a holomorphic chart $z_1,\cdots,z_n$ so that $X \setminus A$ is equal to
$z_1z_2 \cdots z_k = 0$. The line bundle $L$ is trivial over this chart and
$s = z_1^{a_1} \cdots z_k^{a_k}$ with respect to this trivialization for some $a_1,\cdots,a_k > 0$.
The metric $\|\cdot\|$ is equal to $e^\rho |\cdot|$ where $|\cdot|$ is the standard Euclidean
metric on our trivialization of $L$.
So $f_A = -\log{\|s\|} = -\rho - \sum_{i=1}^k a_i \log{|z_i|}$ in the above coordinate chart.
We have that $\theta_A = -d^c f_A$ and $\omega_A = -dd^c f_A$.
So $df_A(X_{\theta_A}) = \|d f_A\|_J^2$ where $\|\cdot\|_J$ is the metric on the real cotangent
bundle of $X$. There is some constant $\gamma>0$ so that $\|df_A\|^2_J \geq \gamma |df_A|^2$
where $|\cdot|$ (by abuse of notation) is the standard Euclidean metric with respect to our coordinate
chart $z_1, \cdots, z_n$. This implies that:
\[df_A(X_{\theta_A}) \geq \gamma \left( -|d\rho|^2 + \sum_{i=1}^k a_k^2 \left|\frac{1}{z_i}\right|^2 \right).\]
We can assume that the functions $|z_i|$ are bounded above by some constant.
This implies that for any $a>0$, $b \in \R$, there is a constant $\kappa > 0$ such that if
$\text{min}(|z_i|,|z_j|) < \kappa$ 
then $|\frac{1}{z_i}|^2+|\frac{1}{z_j}|^2 \geq a{\log(|z_i|)\log(|z_j|)}+b$.
Because $f_A = -\rho - \sum_{i=1}^k a_i \log{|z_i|}$ we have by the previous two inequalities that
\begin{equation} \label{eqn:differentialinequality}
 df_A(X_{\theta_A}) \geq f_A^2
\end{equation}
sufficiently near $D$.
Now lets look at a flowline $x(t)$ of $X_{\theta_A}$ near $D$.
We have that $y(t) := f_A(x(t))$ satisfies
\[ \frac{dy}{dt} \geq y^2\]
by equation (\ref{eqn:differentialinequality}).
Solving such a differential inequality gives us:
\[ y \geq \frac{y(0)}{1-y(0)t} \]
whose solution blows up in time less than $\frac{1}{y(0)}$ (as we can assume $y(0) > 0$ because we are near $D$).
This implies that if we are inside the region $f_A^{-1}(1,\infty)$ and also sufficiently near $D$
then every flowline of $X_{\theta_A}$ flows off to infinity in time less than $1$.
Hence $(A,\theta_A)$ is a strongly bounded finite type Liouville manifold.
So by Lemma \ref{lemma:boundedembedding} we have an exact symplectic embedding of $(A,\theta_A)$
into $(M,\theta) := (f_A^{-1}(-\infty,C],\nu \theta_A)$ where $C,\nu \gg 1$.
This embedding is a homotopy equivalence.

By Lemma \ref{lemma:affineLiouvilledomain},
$\overline{A}$ is Liouville deformation equivalent to $f_A^{-1}(-\infty,C]$ for any $C \gg 0$
which in turn is Liouville deformation equivalent to $(M,\theta)$.
\qed

The following lemma is technical and will be useful here when $J$
is standard and also useful later on.

\begin{lemma} \label{lemma:connecteddomain}
Let $A$ be a smooth affine variety and $X$ a smooth projective variety compactifying $A$.
We let $J$ be an almost complex structure on $X$ which agrees with the standard complex structure
on $X$ near $X \setminus A$.
Let $u : S \rightarrow X$ be a $J$ holomorphic map where $S$ is a compact nodal Riemann  surface such that no component
of $S$ maps entirely in to $X \setminus A$.
Let $\phi_{s,\|\cdot\|}$ be some plurisubharmonic function associated to an ample line bundle $L$ on $X$.
Then near $u^{-1}(X \setminus A)$ we have that $\phi_{s,\|\cdot\|} \circ u$ has no singularities.
\end{lemma}
\proof of Lemma \ref{lemma:connecteddomain}.
We will define $\phi := \phi_{s,\|\cdot\|} \circ u$.
Near $u^{-1}(X \setminus A)$ there is a holomorphic section
$s$ of $u^* L$ whose zero set is $u^{-1}(X \setminus A)$ such that $\phi = -\log{\|s\|}$.
Let $p \in u^{-1}(X \setminus A)$. We want to show that $\phi$ has no singularities near $p$.
After trivializing $u^* L$ near $p$ we have that $s = g(z) z^l$ where $z$ is some coordinate function
on $S$ near $p = \{z=0\}$ and $g$ is a non-zero holomorphic function near $p$.
Also $\|s\| = e^{-\psi} |s|$ where $|\cdot|$ is the standard Euclidean metric on the trivialization of $L$
and $\psi$ is some function.
If we choose polar coordinates $z = re^{i\vartheta}$ then $\phi = \psi - l\log{r} - \log{|g(z)|}$.
So $-\frac{\partial}{\partial r}(\phi)$ tends to infinity as $r$ tends to $0$
because $-\frac{\partial}{\partial r}(\psi - \log{|g(z)|})$ is bounded but
$\frac{\partial}{\partial r}(l\log{r}) = \frac{l}{r}$.
Hence $\phi$ has no singularities near $u^{-1}(X \setminus A)$.
\qed

The next lemma shows us how to relate holomorphic curves inside smooth affine varieties
with algebraic curves inside compactifications of these smooth affine varieties.
This technique is a degeneration technique.

\begin{lemma} \label{lemma:compactnessresultforvarieties}
Suppose we have a morphism $\pi : Q \rightarrow \C$ of smooth quasiprojective varieties
with the following properties:
\begin{enumerate}
\item There is a symplectic form on $Q$ compatible with the complex structure.
\item The central fiber $\pi^{-1}(0)$  is equal to $F \cup E$ where $F$ and $E$ have the same dimension
and where $F$ is a projective variety.
\item $Q  \setminus E$ is isomorphic to a product $B \times \C$ where $B$ is a smooth affine variety.
The morphism $\pi$ under the above isomorphism is the projection map $B \times \C \twoheadrightarrow \C$.
\item There is a sequence $x_i \in \C \setminus \{0\}$ tending to zero as $i$ tends to infinity
and a holomorphic map $u_{x_i} : S_{x_i} \rightarrow \pi^{-1}(x_i)$ for each $i$.
Here $S_{x_i}$ is a smooth genus zero Riemann surface with $|H_1(S_{x_i},\Q)| \leq k-1$.
This map is not necessarily a proper map and it has energy bounded above by some constant $\Lambda$ with respect to $\omega$
where $\Lambda$ is independent of $i$.
\item There is a neighbourhood $N$ of $F$ with the property that
$u_{x_i}|_{u_{x_i}^{-1}(N)}$ is a proper map for $i$ sufficiently large.
\item All these curves $u_{x_i}$ pass through some point $p_{x_i} \in \pi^{-1}(x_i)$
where $p_{x_i}$ tends to some point $p \in F \setminus E$ as $i$ tends to $\infty$.
\end{enumerate}
Then there is a non-trivial holomorphic curve $v : \P^1 \rightarrow F$
with the property that $v^{-1}(E)$ is at most $k$ points
in $\P^1$. Also $p$ is contained in the image of $v$.
If all the curves $u_{x_i}$ are contained inside some closed subvariety $V$ of $Q$
then $v$ is also contained inside $V$.
\end{lemma}
\proof of Lemma \ref{lemma:compactnessresultforvarieties}.
Choose real compact codimension zero submanifolds with boundary $N_1,N_2$ of $N$ with the property that
the interior of $N_1$ contains $F$ and the interior of $N_2$ contains $N_1$.
We can perturb the boundaries of these manifolds $N_1$ and $N_2$ so that they are transverse to $u_{x_i}$
for all $i$.
Consider the holomorphic curves $u'_{x_i} := u_{x_i}|_{u_x^{-1}(N_2)}$.
By a compactness argument \cite{fish:compactness} using the manifold $N_2$ and curves $u'_{x_i}$,
we have (after passing to a subsequence),
a sequence of compact subcurves $\widetilde{S}_i \subset S_{x_i}$ with the following properties:
\begin{enumerate}
\item the boundary of $\widetilde{S}_i$ is sent by $u_{x_i}$ outside $N_1$.
\item There is a compact surface $\widetilde{S}$ with boundary and a sequence of diffeomorphisms
$a_i : \widetilde{S} \rightarrow \widetilde{S}_i$ so that $u_{x_i} \circ a_i$
$C^0$ converge to some continuous map $v' : \widetilde{S} \rightarrow Q$.
This continuous map is smooth away from some union of curves $\Gamma$ in the interior of $\widetilde{S}$
and $u'_{x_i} \circ a_i$ converge in $C^{\infty}_{\text{loc}}$ to $v'$ outside $\Gamma$.
\item The map $v'$ is equal to $v'' \circ \phi$ where $v''$ is a holomorphic map from
a nodal Riemann surface $S$ with boundary to $N_2$ and $\phi : \widetilde{S} \rightarrow S$ is continuous
surjection, and a diffeomorphism onto its image away from $\Gamma$.
The curves $\Gamma$ get mapped to the nodes of $S$ under $\phi$.
The map $v''$ sends the boundary of $S$ outside $N_1$.
\end{enumerate}
Because each of the curves $u_{x_i}$ are contained in $\pi^{-1}(x_i)$
and $x_i$ tends to zero we also get that the image of $v''$ is contained in $\pi^{-1}(0) = F \cup E$.
Also because the points $p_{x_i}$ converge to $p \in F$, we have some $s \in S$
such that $v''(s) = p$.
Because $v''$ sends the boundary of $S$ outside $N_1$ and that the interior of $N_1$ contains $F$,
we have an irreducible component $\P^1 = K$
inside our nodal Riemann surface $S$ where $v''$ maps $K$ to $F$. Also we can assume that $s \in K$.
Our map $v$ is defined as $v''$ restricted to $K = \P^1$.

We now wish to show that $v^{-1}(E)$ is at most $k$ points.
There is a natural exhausting plurisubharmonic function $\rho$
on $B$ which we can construct using Definition \ref{defn:canonicasymplecticformassociatedtoA}.
We pull back $\rho$ to $\widetilde{\rho}$
under the natural projection map $P_B : B \times \C \twoheadrightarrow B$.
Here we identify $B \times \C$ with $Q \setminus E$.
We know that $v^{-1}(E)$ is the disjoint union of $l$ points for some $l$.
We want to show that $l \leq k$.
The image $v(K)$ of $v$ is a one dimensional projective subvariety of $F$,
and $v(K) \setminus E$ is a (possibly singular) affine subvariety of $B$.
By Lemma \ref{lemma:connecteddomain}
we have that $\rho$ restricted to $v(K)$ is non-singular outside a compact set.
Hence for $C \gg 1$ we have that $\rho^{-1}(C)$ is transverse to $v$
and that $(\rho \circ v)^{-1}(C)$ is a disjoint union $l$ circles.
Also by making $C$ large enough, we have that $v^{-1}(\rho^{-1}(-\infty,C])$ is connected.
Because the maps $u_{x_i} \circ a_i$ converge in $C^\infty_{\text{loc}}$ to $v$
near these $l$ circles for $i$ large enough and that these maps $C^0$ converge,
we get that the connected component $S'_i$ of $(\rho \circ u_{x_i})^{-1}\left((-\infty,C]\right)$
passing through $p_x$ has $l$ boundary components for $i \gg 1$.
Because 
\begin{enumerate}
\item each curve $u_{x_i}$ maps to a smooth affine variety $\pi^{-1}(x_i)$ and
\item $(\rho \circ u_{x_i})^{-1}((-\infty,C])$ is compact for $i$ large enough and
\item $|H^1(S_{x_i},\Q)| \leq k-1$,
\end{enumerate}
we have by Lemma \ref{lemma:topologyofcurve} that
$H_1(S'_i,\Q)$ has rank at most $k-1$ which implies that
$S'_i$ has at most $k$ boundary circles for $i \gg 1$.
But we know for $i$ sufficiently large that it is also has $l$ boundary components, hence $l \leq k$.
This implies that $v^{-1}(E)$ is a union of $\leq k$ points and passes through $p$.

Now suppose that all of these curves $u_{x_i}$ are contained inside some closed subvariety $V$.
Because $u_{x_i} \circ a_i$ $C^0$ converges to $v'$ we have that $v'$ is a subset of $V$
because $V$ is a closed subset of $Q$.
The image of $v$ is contained inside the image of $v'$ which implies that the image of $v$
is contained inside $V$.
\qed

\proof of Theorem \ref{theorem:kuniruledimpliesalgebraicallyuniruled}.
Let $X'$ be some compactification of $A$ by a projective variety.
By the Hironaka resolution of singularities theorem \cite{hironaka:resolution}
we can resolve the compactification $X'$ of $A$
so that it is some smooth projective variety $X$ with the property that $X \setminus A$
is a smooth normal crossing divisor.
This variety can be embedded into $\P^N$ so that $X \setminus A$ is equal
to $X$ intersected with some linear hypersurface $\P^{N-1}$ in $\P^N$.
So we can view $A$ as a subvariety of $\C^N = \P^N \setminus \P^{N-1}$.
We will let $D_X$ be the effective ample divisor given by restricting
$\P^{N-1}$ to $X$.

We start with $\P^1 \times \P^N$. The divisor
$D := \{\infty\} \times \P^N + \P^1 \times \P^{N-1}$ is ample.
Let $P := \text{Bl}_{\bracket{0} \times \P^{N-1}} \P^1 \times \P^N$
be the natural blowup map along $\{0\} \times \P^{N-1}$ and let
$\widetilde{D}$ be the proper transform of $D$ in $P$. We let $E$ be the exceptional divisor.
Then $k \widetilde{D} + (k-1)E$ is ample inside $P$ for $k \gg 1$.
Let $\pi : P \rightarrow \P^1$ be the composition of the blowdown map with the projection map to $\P^1$.
The fiber $\pi^{-1}(0)$ is a union of two divisors $F + E$ and this is linearly equivalent
to $\pi^{-1}(\infty)$. Hence $E$ is linearly equivalent to $\pi^{-1}(\infty)-F$.
Let $D'$ be the divisor $k\widetilde{D} + (k-1)(\pi^{-1}(\infty) - F)$.
This is ample and the associated line bundle $L_{D'}$ admits a metric $\|\cdot\|$ whose curvature form
is a positive $(1,1)$-form. This gives us a symplectic form on $X$.
Let $s$ be a meromorphic section of $L_{D'}$ so that $s^{-1}(0) - s^{-1}(\infty) = D'$.
We have that $-d^c \log{\|s\|}$ restricted to $\pi^{-1}(x) \setminus \text{support}(D')$ ($x \neq 0$)
makes this fiber into a Liouville manifold.
We have that $D'$ is the disjoint union of $D'_1 := k\widetilde{D} + (k-1)\pi^{-1}(\infty)$
and $-(k-1)F$.
Also $-\log{\|s\|}$ tends to $+\infty$ as we approach $D'_1$
and $-\infty$ as we approach $F$.
Hence
$P_C := \left( (-\log{\|s\|})^{-1}((-\infty,C]) \right) \cup F$ is a compact submanifold
of $X \setminus \text{support}(D'_1)$ for generic $C \gg 1$ whose interior contains $F$.

Consider $\P^1 \times X \subset \P^1 \times \P^N$ and let $P_X$ be the proper
transform of $\P^1 \times X$ inside $P$. 
We let $\pi_X$ be the restriction of $\pi$ to $P_X$.
We have $A_x := \pi_X^{-1}(x) \setminus \text{support}(D'_1)$
are all isomorphic smooth affine varieties when $x \neq 0$.
Also if $X_x := \pi_X^{-1}(x)$ then these isomorphisms extend to isomorphisms
$\phi_{x,y} : X_x \rightarrow X_y$ so that
$\phi^* L|_{X_y} = L|_{X_x}$. All these affine varieties are isomorphic to $A$.
So by Lemma \ref{lemma:algebraicvarietysymplecticforminvariance} we have that
all these affine varieties are symplectomorphic with respect to the symplectic form
$-dd^c \log{\|s\|}$.
Combining this with Lemma \ref{lemma:affinevarietycontainedinliouvilledomain}
we then get that all these varieties can be codimension $0$ symplectically embedded into
a fixed Liouville domain $(M,\theta)$ which is Liouville deformation equivalent to $\overline{A}$.
Also these embeddings are homotopy equivalences.
Because $\overline{A}$ is $(k,\Lambda)$ uniruled, we have by Lemma \ref{corollary:deformationinvariance}
that $(M,\theta)$ is $(k,\Lambda')$ uniruled for some $\Lambda'>0$.
We define $P_A \subset P_X$ to be equal to $P_X \setminus (\text{support}(D'_1) \cup E)$.
This is isomorphic to $\C \times A$.
Let $p$ be any point in $P_X \cap (F \setminus E)$ and let $x_i \in \C \setminus \{0\}$,
$p_{x_i} \in P_C \cap A_{x_i}$ 
be a family of points in $A_{x_i}$ which all converge to $p$ as $i$ tends to $\infty$.
For every $x_i$ choose a Liouville domain $N_{x_i}$ which is an exact codimension $0$
symplectic submanifold of $A_{x_i}$
containing $P_C \cap A_{x_i}$ and so that the embedding map $N_{x_i} \hookrightarrow A_{x_i}$ is a homotopy equivalence.
By Lemma \ref{theorem:subdomainuniruled}, we have that $N_{x_i}$ is $(k,\Lambda')$ uniruled
because it can be symplectically embedded into $M$ so that the embedding is a homotopy equivalence.
So for each $i$ there is a proper $J$ holomorphic curve $u_{x_i} : S_{x_i} \rightarrow N_{x_i}^0$
where $u_{x_i}$ has energy $\leq \Lambda'$ and $|H_1(S_{x_i},\Q)|\leq k-1$.
In particular $u_{x_i}|_{u_x^{-1}(P_C^o)}$ are properly embedded holomorphic curves inside
the interior $P_C^o$ of the compact manifold $P_C$.
By Lemma \ref{lemma:compactnessresultforvarieties},
there is an algebraic map $v : \P^1 \rightarrow P_X \cap F$ with the property that
$v(q) = p$ for some $q \in \P^1$ and $v^{-1}(E)$ is a union of at most $k$ points.
After identifying $A$ with $P_X \cap (F \setminus E)$ we then get that  $v|_{v^{-1}(A)}$
is also algebraic and passing through $p$ and $v^{-1}(A)$
is $\P^1$ with at most $k$ punctures.
Hence $A$ is algebraically $k$-uniruled.
\qed

\section{Introduction to Gromov Witten invariants} \label{section:GWintro}

Genus $0$ Gromov Witten invariants for general symplectic manifolds
have now been defined in many different ways:
\cite{FukayaOno:Arnold}, \cite{CieliebakMohnke:symplectichypersurfaces},
\cite{HoferWysockiZehnder:polyfoldapplications1} and \cite{LiTian:sympGW}.
Earlier work for special symplectic manifolds such as projective varieties
of complex dimension $3$ or less are done in
\cite{Ruan:topologicalsigma}, \cite{Ruan:3folds} and
\cite{RuanTian:quantum}.
Many of the applications of this paper appear in complex dimension $3$
or less.
These invariants can also be defined in a purely algebraic way
 \cite{LiTian:algGW}, \cite{BehrendFantechi:normalcone}
and \cite{Behrend:GW} but we will not use these theories here.
We will use the Gromov Witten invariants defined for general symplectic manifolds.
All of our calculations are done for complex structures where all of the curves
in the relevant homology class are regular and unobstructed
(and also somewhere injective) and so are relatively easy calculations.
Also most of these (or similar) calculations have been done before
in \cite{Mcduff:rationalruled}, \cite{Ruan:virtual} and \cite{Kollar:lowdegree}.

We start with a compact symplectic manifold $X$, a natural number $k$ and an element $\beta \in H_2(X)$.
Let $d := 2 \left( n - 3 + k + c_1(X).\beta \right)$ where $n$ is half the dimension of $X$.
Choose $k$ cohomology classes $\alpha_1,\cdots,\alpha_k \in H^*(X,\Q)$ so that the sum of their degrees is $d$.
For any compatible almost complex structure one has the set
${\mathcal M}(\beta,J,k)$ of $J$ holomorphic maps $u : S \rightarrow X$
where $S$ is a genus $0$ compact nodal Riemann surface with $k$ labeled marked points.
This nodal curve has to be stable which means that if an irreducible component of this surface
maps to a point then that component must have at least three of these marked points.
There are natural maps $\text{ev}_i : {\mathcal M}(\beta,J,k) \to X$
which send a curve $u : S \to X$ to $u(x_i)$ where $x_i$ is the $i$th marked point in $S$.
It turns out (in nice circumstances) that ${\mathcal M}(\beta,J,k)$
is a topological space with a homology class
\[\left[ {\mathcal M}(\beta,J,k) \right]^{\text{vir}} \in H_d({\mathcal M}(\beta,J,k),\Q).\]
One then has
\[ \langle \alpha_1, \cdots, \alpha_k \rangle^X_{0,\beta} :=
\int_{\left[ {\mathcal M}(\beta,J,k) \right]^{\text{vir}} }
\text{ev}_1^* \alpha_1 \wedge \cdots \wedge \text{ev}_k^* \alpha_k.\]

The genus $0$ Gromov Witten invariant
\[ \langle \alpha_1, \cdots, \alpha_k \rangle^X_{0,\beta} \in \Q\]
satisfies the following properties:
\begin{enumerate}
\item \label{item:existenceofcurves}
If $\langle \alpha_1, \cdots, \alpha_k \rangle^X_{0,\beta} \neq 0$
for some $\alpha_1,\cdots,\alpha_k$ then for every compatible $J$, there exists
a $J$ holomorphic map $u : S \rightarrow X$ from a genus $0$ nodal curve $S$
representing the class $\beta$.
\item \label{item:smoothcalculation}
Suppose that $X$ is a smooth projective variety with its natural complex structure $J$.
Suppose that every rational curve $C$ representing the class $\beta$ is smooth, embedded,
and satisfies $H^1(C,T_X|_C) = 0$ where $T_X$ is the tangent sheaf.
Then $\langle \alpha_1, \cdots, \alpha_k \rangle^X_{0,\beta} \neq 0$
for some $\alpha_1,\cdots,\alpha_k$.
\end{enumerate}
The reason why (\ref{item:smoothcalculation}) is true is that
${\mathcal M}(\beta,J,k)$ in this case is a complex manifold
of dimension $d$ for every $k$ and $\left[ {\mathcal M}(\beta,J,k) \right]^{\text{vir}}$
is equal to its fundamental class.
For $k$ large enough, the map:
\[ \text{ev}_1 \times \cdots \times \text{ev}_k : {\mathcal M}(\beta,J,k) \to X^k \]
is a holomorphic map which is a branched cover onto its image.
If we restrict the natural product symplectic structure
$\omega_{X^k} := \omega_1 + \cdots + \omega_k$ on $X^k$
to ${\mathcal M}(\beta,J,k)$ then it is also a symplectic structure on this moduli space
away from the branching locus.
Hence $\omega^d_{X^k}$ restricted to ${\mathcal M}(\beta,J,k)$
is a positive multiple of the volume form on an open dense subset and so
it evaluates non-trivially with the fundamental class.
In particular we have that $\omega_1^{i_1} \wedge \cdots \wedge \omega_k^{i_k}$
evaluates non-trivially with the fundamental class for some $i_1,\cdots,i_k$.
So if we choose $\alpha_l := \omega_1^{i_l}$ then
$\langle \alpha_1, \cdots, \alpha_k \rangle^X_{0,\beta} \neq 0$.
This argument is almost exactly the same as an argument at the end of the proof of
\cite[Theorem 4.10]{RuanJianxunTianJun:birational}.

\section{Uniruledness criteria for affine varieties} \label{section:affineuniruled2}

In this section we will give another definition of uniruledness for smooth affine varieties.
The main theorem of this section is to show that any smooth affine variety satisfying this uniruledness
condition has an associated Liouville domain which is also $(k,\Lambda)$-uniruled.

We say that a smooth affine variety $A$ is {\bf compactified $k$-uniruled} if $A$
has some compactification $X$ by a smooth projective variety so that if
$D = X \setminus A$ then we have the following properties:
\begin{enumerate}
\item There is an effective ample divisor $D_X$ on $X$
whose support is $D$ and whose associated line bundle has a metric whose curvature form gives
us some symplectic form $\omega_X$ on $X$.
\item Let $J$ be an almost complex structure compatible with $\omega_X$ which
is the standard complex structure near $D$.
Then there is a dense set $U_J \subset A$ so that for every point $p \in U_J$
such that $J$ is integrable near $p$ we have a $J$ holomorphic map
$u : \P^1 \rightarrow X$ passing through $p$ such that $u^{-1}(D \setminus A)$ is a union
of at most $k$ points.
\item The energy of this curve is bounded above by some fixed constant $\Lambda$.
\end{enumerate}

We will now give some easier criteria for being compactified $k$-uniruled.
Our symplectic form on $X$ comes from some ample divisor.
\begin{lemma} \label{lemma:fibrationgwkuniruled}
Suppose that we have a morphism $\pi : X \rightarrow B$ whose generic fiber is $\P^1$
where the base $B$ is projective.
Let $\beta \in H_2(X)$ be the class of this curve.
Then for every compatible almost complex structure $J$
which is integrable on some open set $U$ containing
a point $p$, there is some $J$ holomorphic curve $u : S \rightarrow X$ passing through $p$
representing the class of the fiber.
Here $S$ is a genus $0$ nodal curve.
\end{lemma}
\proof of Lemma \ref{lemma:fibrationgwkuniruled}.
Let $F$ be any regular fiber of $\pi$. This is isomorphic to $\P^1$.
Blow up $X$ to $\widetilde{X}$ at some point in $F$.
Let $\widetilde{F}$ be the proper transform of $F$ inside $X$
and  let $\widetilde{\beta} \in H_2(\widetilde{X},\Q)$ be its respective homology class.
The only curve in this homology class is $\widetilde{F}$.
If we restrict the tangent bundle $T_{\widetilde{X}}$ to this curve then it is isomorphic to:
${\mathcal O}(2) \oplus {\mathcal O}(-1)^{\oplus n-1}$.
Hence $H^1(\widetilde{F},T_{X}|_{\widetilde{F}}) = 0$.
By property (\ref{item:smoothcalculation}) there exists
$\alpha_1,\cdots,\alpha_k$ such that
$\langle \alpha_1, \cdots, \alpha_k \rangle^X_{0,\widetilde{\beta}} \neq 0$.

Let $p_i$ be a sequence of points in $U$ converging to $p$ where $p_i$
is contained inside a smooth fiber $F_i$.
Because $J$ is integrable in $U$, we can blowup $X$ at $p_i$
giving us a new symplectic manifold $X_i$ along with a compatible almost complex structure
so that the blowdown map
is holomorphic. Let $\widetilde{F}_i$ be the proper transform of $F_i$ in $X_i$
and $\beta_i \in H_2(X_i,\Z)$ its respective homology class.
Then because $\langle \alpha_1, \cdots, \alpha_k \rangle^X_{0,\beta_i} \neq 0$ for some cohomology
classes $\alpha_i$,
we have  by property (\ref{item:existenceofcurves}) a $J$ holomorphic curve $u'_i : S_i \rightarrow X_i$
representing $\beta_i$. By composing this map with the blowdown map,
we get a $J$ holomorphic curve $u_i : S_i \rightarrow X$ passing through $p_i$
and representing $\beta$.
By a Gromov compactness argument one then gets a holomorphic curve $u : S_i \to X$
passing through $p$ representing the class $\beta$.
\qed

\begin{lemma} \label{lemma:nefdivisoruniruled}
Suppose that we have a morphism $\pi : X \rightarrow B$ whose generic fiber
is $\P^1$ where $B$ is a projective variety.
Suppose also that $D'$ is an effective nef divisor whose support is equal to $D$
with the property that $\beta.D' \leq k$. Then $A$ is compactified $k$-uniruled.
\end{lemma}
\proof of Lemma \ref{lemma:nefdivisoruniruled}.
Choose any effective ample divisor $D_X$ whose support is $D$
and let $\omega_X$ be the symplectic form associated to this divisor.
Let $J$ be any compatible almost complex structure which is standard near $D$.
Let $p$ be any point in $A$ where $J$ is integrable near $p$. 
By Lemma \ref{lemma:fibrationgwkuniruled}, there is a $J$ holomorphic curve
$u : S \rightarrow X$ representing the homology class of the fiber passing through $p$.
Here $S$ is a nodal curve with irreducible components $S_1,\cdots,S_l$.
Let $S_i$ be any component passing through $p$.
Then by positivity of intersection we have that $u(S_i).D \leq u(S_i).D'$
and $u(S_i).D' \leq \sum_i u(S_i).D' = u(S).D' \leq k$.
Hence $u(S_i).D \leq k$.
The energy of the curve $u|_{S_i}$ is bounded above by $\beta.D_X$.
Because $S_i$ is isomorphic to $\P^1$ we then get that $A$ is compactified $k$-uniruled.
\qed

\begin{theorem} \label{theorem:algebraicuniruledimplesuniruled}
Suppose that $A$ is a smooth affine variety that is compactified $k$-uniruled.
Let $\overline{A}$ be its associated Liouville domain. Then $\overline{A}$ is
$(k,\Lambda)$-uniruled for some $\Lambda$.
\end{theorem}

Before we prove this theorem we need a lemma and a definition.
Let $X$ be a smooth projective variety with a smooth normal crossing divisor
$D$ so that $X \setminus D$ is affine.
A map $u : S \to X$ is said to be a {\it $k$-curve} if every irreducible
component $\Sigma$ of $S$ either maps to $D$, or $u^{-1}(D) \cap \Sigma$
is a finite set of size at most $k$.


\begin{lemma} \label{lemma:compactnessresult}
Let $A$ be a smooth affine variety and $X$ a smooth projective variety compactifying $A$.
We equip $X$ with a symplectic form $\omega_{\|\cdot\|}$ coming from some ample line bundle.
We let $J$ be a compatible almost complex structure on $X$ which agrees with the standard complex structure
on $X$ near $X \setminus A$.
Let $u_i : S_i \rightarrow X$ be a sequence of $J$ holomorphic maps
where $S_i$ is a connected genus $0$ nodal Riemann surface and
where all the $u_i$'s have energy bounded above by some fixed constant.
If the $u_i$ are all $k$ curves and Gromov converge to $u : S \rightarrow X$
then $u$ is also a $k$ curve.
\end{lemma}
\proof of Lemma \ref{lemma:compactnessresult}.
Let $S_i^o := u_i^{-1}(A)$ and $S^o := u^{-1}(A)$.
We want to show that the rank of $H_1$ of each connected component of $S^o$ is at most $k-1$.
Because the $u_i$ Gromov converge, that means that there is a smooth real surface $\widetilde{S}$
and a series of continuous maps $\alpha_i : \widetilde{S} \rightarrow S_i$, $\alpha : \widetilde{S} \rightarrow S$
satisfying:
\begin{enumerate}
\item $\alpha_i$ and $\alpha$ are diffeomorphisms away from a $1$-dimensional submanifold $\Gamma \subset \widetilde{S}$
and away from the nodes of $S_i$ and $S$.
\item  $\Gamma$ maps to the nodes of $S$ under $\alpha_0$.
\item $u_i \circ \alpha_i$ $C^0$ converge to $\alpha_0 \circ u$ and these maps $C^\infty_{\text{loc}}$ converge
away from $\Gamma$.
\end{enumerate}
Choose an exhausting plurisubharmonic function $\phi : A \rightarrow \R$ associated to some line bundle
$L$, section $s$ and metric $\|\cdot\|$ on $L$.
Gromov convergence means that for $c \gg 1$ and $i \gg 1$
we have that every node of $S_i^o$ and also $\alpha_i(\Gamma) \cap S_i^o$
is mapped via $u_i$ to $\phi^{-1}(-\infty,c]$.
We can also assume that the same is true for $S^0 := u^{-1}(A)$.
For large enough $i$ and for generic $c$ large enough we have
$u_i$ is smooth near $\phi^{-1}(c)$ and also transverse to this hypersurface.
We can assume the same properties hold for $u$.
Let $\Sigma_i$ be a sequence of connected components of $S_i^o$
which converge to a connected component $\Sigma$ of $S^o$.
We have that $\Sigma_i \cap u_i^{-1}(\phi^{-1}(c))$ is a union of $l_i$ smooth circles in $S_i$
and $\Sigma \cap u^{-1}(\phi^{-1}(c))$ is a union of $l$ circles for some $l_i,l$.
This means that $H_1(\Sigma_i \cap u_i^{-1}(\phi^{-1}(-\infty,c))$ has rank $\leq l_i-1$.
By Lemma \ref{lemma:topologyofcurve} we then get that $H_1(\Sigma \cap u_i^{-1}(\phi^{-1}(-\infty,c))$ has rank
less than or equal to $|H_1(\Sigma_i)| \leq k-1$.
Hence $l_i \leq k-1$ for all $i$. So $\Sigma_i \cap u_i^{-1}(\phi^{-1}(c))$ is a union of at most $k$ circles.
Because $u_i \circ \alpha_i$ $C^\infty$ converge to $u_i \circ \alpha$ near $\phi^{-1}(c)$
we get that $\Sigma \cap u^{-1}(\phi^{-1}(c))$ is also a union of at most $k$ circles.
This is true for all $c$ sufficiently large.
Hence $\text{rank}(H_1(\Sigma))\leq k-1$ for each connected component $\Sigma$ of $S^o$.
Hence $u$ is a $k$ curve.
\qed

\proof of Theorem \ref{theorem:algebraicuniruledimplesuniruled}.
Because $A$ is compactified $k$-uniruled, we have a compactification $X$ with divisor $D$
so that:
\begin{enumerate}
\item There is an effective ample divisor $D_X$ on $X$
with support $D$ whose associated line bundle has a metric with curvature form
$\omega_X$ on $X$. Here $\omega_X$ is a symplectic form.
\item Let $J$ be an almost complex structure compatible with $\omega_X$ which
is the standard complex structure near $D$.
Then there is a dense set $U_J \subset A$ so that for every point $p \in U_J$
such that $J$ is integrable near $p$ we have a $J$ holomorphic map
$u : \P^1 \rightarrow X$ passing through $p$ which is a $k$ curve.
\item The energy of this curve is bounded above by some fixed constant $\Lambda'$.
\end{enumerate}
We have a plurisubharmonic function $\rho := -\log{|s|}$ where $s$ is a section of $L$ with
$s^{-1}(0) = D_X$.
For $c \gg 1$ we have that $A_c := \rho^{-1}(-\infty,c]$ is a Liouville domain deformation equivalent to
$\overline{A}$ by Theorem \ref{lemma:affineLiouvilledomain}.
We now let $J$ be any almost complex structure which coincides with the standard one near $D$
and coincides with any convex almost complex structure inside $A_c$.
Let $p$ be any point in the interior of $A_c$ where $J$ is integrable on a neighbourhood of $p$.
Choose a sequence of points $p_i \in U_J$ converging to $p$.
There is a map $u_i : \P^1 \to X$ of energy bounded above by $\Lambda'$ passing through $p_i$
so that $u_i$ is a $k$ curve.
There is a subsequence which Gromov converges to a map $v : S \rightarrow X$ of energy bounded above by $\Lambda'$
passing through $p$.
Here $S$ is a genus $0$ nodal curve. 
By Lemma \ref{lemma:compactnessresult}, we then get that $v$ is a $k$ curve.
Let $S'$ be an irreducible component of $S$ whose image under $v$ contains $p$,
$S'' := S' \cap v^{-1}(A)$ and $\Sigma := S'' \cap v^{-1}(A_c^0)$ where
$A_c^0$ is the interior of $A_c$.
By Lemma \ref{lemma:topologyofcurve} we have that $|H_1(\Sigma,\Q)| \leq k$
because $|H_1(S'',\Q)| \leq k$. Let $u := v|_\Sigma$.
The energy of $u$ is bounded above by $\Lambda'$.
This implies that $A_c$ is $(k,\Lambda')$-uniruled.
By Corollary \ref{corollary:deformationinvariance} we then get that $\overline{A}$
is $(k,\Lambda)$-uniruled for some $\Lambda>0$.
\qed

\section{Log Kodaira dimension and uniruledness} \label{section:logkodaira}

We will now define log Kodaira dimension.
Let $L$ be any line bundle on a projective variety $X$.
If $L^{\otimes k}$ has no global sections for any $k$ then we define $\kappa(L) := -\infty$.
Otherwise $L^{\otimes k}$ defines a rational map from $X$ to $\P(H^0(L^{\otimes k}))$ for some $k$.
We define $\kappa(L)$ in this case to be maximum  dimension of the image of this map
over all $k$ where this map is defined.
The number $\kappa(L)$ is called the {\it Kodaira dimension of} $L$.
If $Q$ is any smooth quasiprojective variety then we define its {\it log Kodaira dimension}
$\overline{\kappa}(Q)$ as follows:
Choose some compactification of $Q$ by a smooth projective variety $X$ so that
the associated compactification divisor $D$ is smooth normal crossing.
The {\it log Kodaira dimension} of $Q$ is defined to be $\kappa(K_X + Q)$
where $K_X$ is the canonical bundle of $X$.
This an invariant of $Q$ up to algebraic isomorphism.

Before we look at smooth affine varieties in dimension $2$ and $3$ we need a lemma relating
uniruledness with log Kodaira dimension.
\begin{lemma} \label{lemma:logkodairaunirulednessrelation}
Suppose that $A$ is algebraically $k$-uniruled. If $k=1$ then
$A$ has log Kodaira dimension $-\infty$ and if $k=2$ then
$A$ has log Kodaira dimension $\leq \text{dim}_\C A - 1$.
\end{lemma}
\proof of Lemma \ref{lemma:logkodairaunirulednessrelation}.
First of all, we compactify $A$ to some smooth projective variety $X$.
Let $D$ be the compactification divisor.
Because $A$ is algebraically $k$-uniruled, we have that $X$ is uniruled by $\P^1$'s.
By using the theory of Hilbert schemes (see \cite{Kollar:rationalcurves}) there is a surjective morphism:
\[ \text{ev} : M \times \P^1 \twoheadrightarrow X \]
where $M$ is a reduced projective variety.
We define  $D_M := \text{ev}^{-1}(D)$.
We let $V$ be the subvariety of $M$ with the property that $q \in M$
is contained in $V$ if and only if
$D_M \cap (\{q\} \times \P^1)$ is a set of size at most $k$.
Because $A$ is $k$ uniruled we can assume that $M$ satisfies: $\text{ev}(V \times \P^1)$ is dense in $X$ .
Hence we have a dominant morphism
$(V \times \P^1) \setminus D_M \twoheadrightarrow A$.
We will define $W$ be equal to $(V^{\text{sm}} \times \P^1) \setminus D_M$
where $V^{\text{sm}}$ is the smooth part of $V$ which is a non-empty Zariski open subset of $V$.
In particular we have a morphism $\pi_W$ from $W$ to $A$ whose image contains a dense open set.
We can choose $V' \subset V$ to be a subvariety of complex dimension $\text{dim}(X) -1$
so that the image $\pi_W((V' \times \P^1) \setminus D_M)$ still contains a dense open subset of $A$.
We define $W'$ to be $(V' \times \P^1) \setminus D_M$.
So $\pi_{W'} := \pi_W|_{W'}$ is a dominant morphism from $W'$ to $A$.
The projection map $W' \twoheadrightarrow V'$ has generic fiber equal to $\P^1$ minus at most $k$ points.
By the Iitaka Easy Addition Theorem (\cite[Theorem 4]{Iitaka:logarithmic},
\cite[Theorem 11.9]{Iitaka:algebraicgeometry}) we have that the log Kodaira dimension
of $W'$ is equal to $-\infty$ if $k=1$ and it is $\leq \text{dim}_\C(A)-1$ if $k = 2$.
Because there is a dominant morphism from $W'$ to $A$,
we have by the logarithmic ramification formula (\cite{Iitaka:logarithmic}, \cite[Theorem 11.3]{Iitaka:algebraicgeometry})
that the log Kodaira dimension of $W'$ is greater than or equal to the log Kodaira dimension of $A$.
Combining the above two facts we have that if $k=1$ then $A$ has log Kodaira dimension $-\infty$
and if $k=2$ then $A$ has log Kodaira dimension $\leq \text{dim}_\C(A)-1$.
\qed

\begin{lemma} \label{lemma:compactifieduniruledlkd}
Suppose that $A$ and $B$ are symplectomorphic smooth affine varieties.
Suppose that $A$ is compactified $k$-uniruled. If $k=1$ then
$B$ has log Kodaira dimension $-\infty$ and if $k=2$ then
$B$ has log Kodaira dimension $\leq \text{dim}_\C A - 1$.
\end{lemma}
\proof of Lemma \ref{lemma:compactifieduniruledlkd}.
By Theorem \ref{theorem:algebraicuniruledimplesuniruled} we get that
the Liouville domain $\overline{A}$ associated to $A$ is $(k,\Lambda)$ uniruled.
Because $B$ is symplectomorphic to $A$ we then get by Lemma \ref{theorem:symplecticinvariance}
that the Liouville domain $\overline{B}$ is $(k,\Lambda')$ uniruled.
So by Theorem \ref{theorem:kuniruledimpliesalgebraicallyuniruled}, $B$ is algebraically $k$ uniruled.
Hence by Lemma \ref{lemma:logkodairaunirulednessrelation}, $B$ has
the required log Kodaira dimension.
\qed

\subsection{Dimension 2}

The aim of this section is to prove:
\begin{theorem} \label{theorem:acyclicsurfaceinvariance}
Let $A$, $B$ be symplectomorphic acyclic smooth affine surfaces.
Then they have the same log Kodaira dimension.
\end{theorem}

Before we prove this we need a compactified uniruled criterion in dimension $2$
and some other preliminary lemmas.

\begin{lemma} \label{lemma:countablymanycurves}
Let $X$ be any compact symplectic manifold of real dimension $4$
and $J$ any almost complex structure compatible with the symplectic form.
Then there is a dense subset of points $U_J \subset X$ with the property
that every $J$ holomorphic map $u : \P^1 \to X$ with $u(\P^1) \cap U_J \neq \varnothing$ satisfies $u_*([\P^1])^2 \geq 0$.
In fact $U_J$ is a countably infinite intersection of open dense subsets.
\end{lemma}
\proof of Lemma \ref{lemma:countablymanycurves}.
Let $E$ be a homology class satisfying $E.E < 0$.
Let $u_i : \P^1 \to X$, $i=1,2$ be two $J$ holomorphic curves representing this class.
We have that $(u_1)_*([\P^1]) \cdot (u_2)_*([\P^1])$ is negative.
By positivity of intersection we then have that the images of $u_1$
and $u_2$ must coincide.
We write $E_J$ for this image.
By Sard's theorem, the complement of this image is open and dense.
The set of images of $J$ holomorphic curves $u : \P^1 \to X$
with negative self intersection number is $\cup_{E \in H_2, E.E < 0} E_J$.
The complement is a countable intersection of open dense subsets,
which is also dense.
\qed

\begin{lemma} \label{lemma:irreducibleholomorphiccurves}
Suppose we have a morphism $\pi : X \rightarrow B$ where $X$
is a smooth projective surface and $B$ is a curve.
Let $\omega_X$ be a symplectic form associated to an effective ample divisor $D_X$
and $J$ a compatible almost complex structure.
Suppose that we have a $J$ holomorphic map $v : \Sigma \to X$
whose fundamental class represents the fiber $[F] \in H_2(X)$,
and with the property that every irreducible component $\Sigma'$ of $\Sigma$
satisfies $v_*([\Sigma']).[F] = 0$.
Then there is a dense set $U_J \subset X$ with the property that every $J$
holomorphic map $u : S \to X$ where $S$ is a connected nodal curve
which intersects $U_J$ and represents $[F]$ has the property that $S$ is
irreducible.
\end{lemma}
\proof of Lemma \ref{lemma:irreducibleholomorphiccurves}.
We choose $U_J$ to be the set of points with the property that
every $J$ holomorphic curve passing through a point in $U_J$
has some irreducible component with non-negative intersection number.
This is dense by Lemma \ref{lemma:countablymanycurves}.
Let $S$ be a union of irreducible components $S_1,\cdots,S_l$.
We will suppose without loss of generality that $S_1 \cap U_J$ is non-empty,
and so that $(u|_{S_1})^2 \geq 0$.
Suppose for a contradiction, $u_*([S_i]).[F] < 0$ for some $i$.
Then by positivity of intersection we have that $u(S_i)$ is contained
inside $v(\Sigma')$ for some irreducible component $\Sigma'$ of $\Sigma$.
Because $v_*([\Sigma']).[F] = 0$, we have
$u_*([S_i]).[F] = 0$ which is a contradiction.
Hence $u_*([S_i]).[F] \geq 0$ for each $i$.
Because $\sum_i u_*([S_i]).[F] = 0$, this implies that $u_*([S_i]).[F] = 0$ for all $i$.

Suppose for a contradiction, $S$ has more than one irreducible component.
Because $S$ is connected we have then that $u_*([S_1]).u_*([S_j]) \neq 0$ for some $j \neq 1$,
and because $(u_*([S_1]))^2 \geq 0$ we then get that $u_*([S_1]).(\sum_i u_*([S_i])) > 0$
which is impossible because $u_*([S]) = [F]$.
Hence $S$ is irreducible.
\qed

\begin{corollary} \label{corollary:uniruledfibrationproperty}
Let $A$ be a smooth affine variety with an SNC compactification $X$,
and let $D$ be the associated compactification divisor.
Suppose we have a morphism $\pi : X \to B$
satisfying the hypotheses of Lemma \ref{lemma:irreducibleholomorphiccurves}
for any compatible $J$ which is standard near $D$.
Suppose that a general fiber of $\pi$ intersects $D$ $k$ times.
Then $A$ is compactified $k$ uniruled.
\end{corollary}
\proof of Corollary \ref{corollary:uniruledfibrationproperty}.
By Lemma \ref{lemma:irreducibleholomorphiccurves},
there is a dense subset $U_J \subset A$ with the property that
any $J$ holomorphic curve $u : S \to X$ representing $F$ passing through $p \in U_J$
has the property that $S$ is irreducible.
Let $p$ be any point in $U_J$ such that $J$ is integrable near $p$, then by Lemma \ref{lemma:fibrationgwkuniruled},
we have that there is such a $J$ holomorphic map passing through $p$.
Because $S$ is irreducible, it intersects $D$ in at most $k$ points.
Putting all of this together gives us that $A$ is compactified $k$ uniruled.
\qed

\begin{lemma} \label{lemma:nefdivisor}
Suppose that $X$ is a smooth projective surface and $D$, $E$ divisors so that:
\begin{enumerate}
\item $E$ is irreducible and $D \cup E$ is smooth normal crossing.
\item If $D'$ is the union of divisors in $D$ not intersecting $E$ then $D'$ is connected
and intersects every irreducible component of $D$.
\item There is an effective nef divisor $G$ whose support is contained in $D'$.
\end{enumerate}
Then there is an effective nef divisor $D_X$ whose support is $D'$ so that:
\begin{enumerate}
\item
For every irreducible curve $C$ in $D$ and not in the support of $G$,
we have $D_X.C > 0$.
\item $D_X.E = 0$.
\end{enumerate}
\end{lemma}
\proof of Lemma \ref{lemma:nefdivisor}.
Suppose that $W$ is any effective nef divisor whose support is in $D'$ and contains $\text{support}(G)$.
Also suppose every irreducible curve $C$ inside $\text{support}(W)$
but not in $\text{support}(G)$ satisfies $C.W > 0$.
Let $C$ be any irreducible curve of $D'$ not contained
in the support of $W$.
Because $D'$ is connected we can assume
that $C.W \neq 0$.
We let $W' := \kappa W + C$ for $\kappa \gg 1$.
This is an effective nef divisor with larger support than $W$.
For $\kappa$ large enough we have $C.W' > 0$.
Hence every irreducible curve $C$ inside $\text{support}(W')$
satisfies $C.W' > 0$ if $C$ is not in $\text{support}(G)$.
Therefore we can construct effective nef divisors starting with $G$
with larger and larger support until we get an effective nef divisor
$D_X$ whose support is equal to $D'$.
Every irreducible curve $C$ in $D'$ intersects $D_X$
positively unless it is in the support of $G$.
Also if $C$ is an irreducible curve in $D$ not contained in $D'$ then it intersects $D'$ and hence $C.D_X>0$.
We also have that $D_X. E = 0$.
\qed

\begin{lemma} \label{lemma:dimension2uniruled3}
Let $A = X \setminus D$ be a smooth affine surface where $X$ is a smooth projective variety
and $D$ is a connected smooth normal crossing divisor.
Suppose that we have a morphism $\pi : X \rightarrow B$ with the following properties:
\begin{enumerate}
\item \label{item:genericratioanlcurve}
The generic fiber is $\P^1$. The base $B$ is a smooth projective curve.
\item \label{item:dotproductbound}
If $F$ is a fiber then $F.D \leq k$ for some $k$.
\item \label{item:specialfibers}
There are two different points $b_1,b_2 \in B$ with the following property:

$\pi^{-1}(b_i) = E_i \cup F_i$ (as reduced curves)
where $E_i$ is an irreducible smooth curve satisfying $E_i.E_i = -1$,
and $F_i$ is reduced.
\item \label{item:normalcrossingcondition}
$F_i \subset D$ but $E_i$ is not contained in $D$.
Also $D \cup E_1 \cup E_2$ is a smooth normal crossing divisor.
\item \label{item:samedotproductcondition}
$E_2.D = E_2.F_2$.
\item \label{item:effectivenefcondition}
There is an effective nef divisor $G$ with the property that
$E_i.G = 0$ and whose support is contained inside $D$.
\item \label{item:connectednesscondition}
The  union of irreducible components of $D$ not containing $E_i$ is connected.
\end{enumerate}
Then $A$ is compactified $k$ uniruled.
\end{lemma}
\proof of Lemma \ref{lemma:dimension2uniruled3}.
Choose any compatible symplectic structure $\omega_X$ coming from an effective
ample divisor $D_X$ whose support is $D$.
Let $J$ be any compatible almost complex structure which is standard near $D$.
We will complete this proof in $3$ steps.
In Step 1 we will construct
$J$ holomorphic curves representing $[E_i]$ such that no irreducible component is contained
in $D$.
In Step 2 we will show that each irreducible component of one the curves from Step 1
has intersection number zero with the fiber.
In Step 3 we will construct our $J$ holomorphic curve passing through $p$
and intersecting $D$ at most $k$ times using Corollary \ref{corollary:uniruledfibrationproperty}.

{\it Step 1}:
By \cite[Lemma 3.1]{Mcduff:rationalruled} there is a $J$ holomorphic map
$u_i : \Sigma_i \to X$ from a connected genus $0$ nodal Riemann surface
$\Sigma_i$ representing the exceptional class $[E_i] \in H_2(X)$.
We assume that no irreducible component of $\Sigma_i$ maps to a point.
Let $\Sigma_i^1,\cdots,\Sigma_i^{l_i}$ be the irreducible components of $\Sigma_i$.
We will now show that $u_i(\Sigma_i^j)$ is not contained in $D$ for each $i,j$.

By using properties (\ref{item:normalcrossingcondition}),(\ref{item:connectednesscondition})
and (\ref{item:effectivenefcondition})
combined with Lemma \ref{lemma:nefdivisor} we have an effective nef divisor $D'_i$ satisfying:
\begin{enumerate}[(a)]
\item \label{item:positiveintersectionwithDX}
For every irreducible curve $C$ in $D$ and not in the support $G$,
we have $D'_i.C > 0$.
\item $D'_i.E_i = 0$.
\end{enumerate}
Suppose for a contradiction that $(u_i)_*([\Sigma_i^y]) \subset D$ for some $y$.
Because $[E_i]^2$ is negative
and the intersection product of $E_i$ with any irreducible component of
$D'_i$ is non-negative we have that $E_i$ cannot be represented
by an effective divisor whose support is in $D'_i$.
If $(u_i)_*([\Sigma_i^y]) \subset D'_i$ then the previous fact tells us
that there is some $\Sigma_i^x$ satisfying $(u_i)_*([\Sigma_i^x]).D'_i \neq 0$.
But this is impossible because $D'_i$ is nef
and $E_i.D'_i = 0$.
Hence $\Sigma_i^y$ is not  contained in $D'_i$, so by property (\ref{item:positiveintersectionwithDX}),
$(u_i)_*([\Sigma_i^y]).D'_i \neq 0$.
This is impossible as $D'_i$ is nef and has intersection number $0$ with $E_i$.
Hence $u_i(\Sigma_i^x)$ is not contained in $D$ for all $i,x$.

{\it Step 2}:
The aim in this step is to show that if $[F]$ is the class of a fiber of $\pi$ then
$(u_1)_*([\Sigma_1^i]).[F] = 0$ for all $i$.
%
Suppose for a contradiction that $(u_1)_*([\Sigma_1^1]).[F] \neq 0$.
Then because $[E_1].[F] = 0$ we have
$(u_1)_*([\Sigma_1^1]).[F] = -\sum_{j=2}^{l_1} (u_1)_*([\Sigma_1^j]).[F]$.
So without loss of generality we can assume that   $(u_1)_*([\Sigma_1^1]).[F] < 0$.
We can represent $[F]$ by $[D_{F_2}] + \kappa (u_2)_*([\Sigma_2])$
by property (\ref{item:specialfibers})
where $D_{F_2}$
is an effective divisor whose support is exactly $F_2$ and $\kappa$ is a positive integer.
Because $(u_1)_*([\Sigma_1^1])$ does not map to $D$, we have by positivity of intersection
that $(u_1)_*([\Sigma_1^1]).[D_{F_2}] \geq 0$.
Hence $(u_1)_*([\Sigma_1^1]). (u_2)_*([\Sigma_2]) < 0$ because $\kappa > 0$.
By positivity of intersection this means that
$u_1(\Sigma_1^1) \subset u_2(\Sigma_2^l)$ for some $l$.
Without loss of generality we will assume that $l=1$.

We have that $E_2.F_2 = E_2.D$ by property (\ref{item:samedotproductcondition})
and that $(u_2)_*([\Sigma_2^i]).[D] \geq 1$ because $A$ is an exact symplectic manifold.
Because $(u_i)_*([\Sigma_i^j])$ is not contained inside $D$ for all $i,j$,
we have $(u_i)_*([\Sigma_i^j]).[F_2] \leq (u_i)_*([\Sigma_i^j]).[D]$.
Using the above two facts, \[(u_2)_*([\Sigma_2^1]).[F_2] =\]
\[\sum_{j=1}^{l_2}(u_2)_*([\Sigma_2^j]).[D] - \sum_{j=2}^{l_2} (u_2)_*([\Sigma_2^j]).[F_2]
\geq (u_2)_*([\Sigma_2^1]).[D] > 0.\]
But this means that $(u_1)_*([\Sigma_1^1]).[F_2] \neq 0$ because
$\varnothing \neq u_1(\Sigma_1^1) \subset u_2(\Sigma_2^1)$.
Hence $E_1.F_2 = (u_1)_*([\Sigma_1]).[F_2] \neq 0$
which is a contradiction because $E_1$ and $F_2$ are in different fibers of $\pi$
by property (\ref{item:specialfibers}).
Hence $(u_1)_*([\Sigma_1^i]).[F] = 0$ for all $i$.

{\it Step 3}:
There is an effective divisor $D_{F_1}$ whose support is $F_1$
and an integer $\kappa'> 0$ with the property that:
$[D_{F_1}] + \kappa' (u_1)_*([\Sigma_1])$ represents $[F]$.
Each irreducible component of the above curve has intersection number zero
with $[F]$ by Step 2 hence by Corollary \ref{corollary:uniruledfibrationproperty},
we get that $A$ is compactified $k$ uniruled.
\qed

\begin{lemma} \label{lemma:dimension2uniruled4}
Let $A = X \setminus D$ be a smooth affine surface where $X$ is a smooth projective variety
and $D$ is a connected smooth normal crossing divisor.
Let $\pi : X \to B$ be a morphism of projective varieties so that the generic fiber is isomorphic to $\P^1$
and intersects $D$ $k$ times.
Let $E$ be a smooth divisor in $X$.
Suppose that:
\begin{enumerate}
\item \label{item:nefdivisornotintersectingE}
There is a nef divisor $G$ whose support is in $D$ such that $E.G=0$.
\item \label{item:divisorintersectingE}
We have $E.E=-1$, $E.D=1$ and $E$ is not contained in $D$.
This means that there is a unique irreducible curve $D_E$ in $D$ intersecting $E$.
We will assume that $D_E.G \neq 0$.
\item \label{item:fiber}
We have that $D_E \cup E$ is contained in a fiber $\pi^{-1}(b)$
and there is an effective divisor $D_F$ whose support is $\pi^{-1}(b) \cap D$
and a natural number $\kappa > 0$ so that $[D_F] + \kappa[E]$ represents the homology class of a fiber of $\pi$.
\end{enumerate}
Then $A$ is compactified $k$ uniruled.
\end{lemma}
\proof of Lemma \ref{lemma:dimension2uniruled4}.
This proof will be done in two steps.
In Step $1$ we will show for any almost complex structure $J$ compatible with the symplectic
form on $X$ and which is standard near $D$, the homology class $[E]$ can be represented by an irreducible $J$
holomorphic curve.  Finally in Step 2 we will use Corollary \ref{corollary:uniruledfibrationproperty}.

{\it Step 1}:
By \cite[Lemma 3.1]{Mcduff:rationalruled} there is a $J$ holomorphic map
$u : \Sigma \to X$ from a connected genus $0$ nodal Riemann surface
$\Sigma$ representing the exceptional class $[E_i] \in H_2(X)$.
Let $\Sigma^1,\cdots,\Sigma^l$ be the irreducible components of $\Sigma$.
In this step we want to show that $l=1$.
Because $u_*([\Sigma^i]).G \geq 0$ for all $i$, and that $u_*([\Sigma]).G = 0$ we then get
$u_*([\Sigma^i]).G = 0$ for all $i$. Hence by property (\ref{item:divisorintersectingE}),
$u(\Sigma^i)$ is not contained in $D_E$ for any $i$.
This means that $u_*([\Sigma^i]).D_E \geq 0$ for all $i$.
The above statement combined with the fact that $D_E.E = 1$ means that
there is exactly one irreducible component $\Sigma^j$ intersecting $D_E$
and this irreducible component intersects $D_E$ with multiplicity $1$.
We may as well assume that $\Sigma^j=\Sigma^1$.
Let $D_X$ be an effective ample divisor whose support is $D$.
Then $E.D_X = u_*([\Sigma]).D_X$, and $u_*([\Sigma^i]).D_X > 0$ for all $i$ because $A$ is an exact symplectic manifold.
Let $c>0$ be the coefficient of $D_E$ in $D_X$.
Then $cE.D_E = E.D_X$ because $D_E$ is the only irreducible divisor in $D$ intersecting $E$.
Also because $u_*([\Sigma^1]).D_E = 1$ we get that $u_*([\Sigma^1]).D_X \geq c$.
Hence $u_*([\Sigma^1]).(D_X - c.D_E) \geq 0$. Also for $i > 1$ we have that
$u_*([\Sigma^i]).D_E = 0$ which implies that $u_*([\Sigma^i]).(D_X - cD_E) = u_*([\Sigma^i]).D_X > 0$.
Hence if $l > 1$ we get that $u_*([\Sigma]).(D_X - cD_E) = \sum_i u_*([\Sigma^i]).(D_X - cD_E) > 0$
which contradicts the fact that $E.(D_X - cD_E) = 0$.
This means that $l = 1$ and so $\Sigma$ is irreducible.

{\it Step 2}:
By Step $1$ we have that $\Sigma$ is irreducible and $u_*([\Sigma])$ represents $E$.
Hence every irreducible component of the $J$ holomorphic curve $D_F \cup u(\Sigma)$
has intersection number $0$ with a fiber of $B$.
So by property (\ref{item:fiber}) combined with Corollary
\ref{corollary:uniruledfibrationproperty} we then get that $A$ is compactified $k$ uniruled.
\qed

\bigskip

We will give a fairly explicit description of acyclic surfaces of log Kodaira dimension
$1$.
The constructions come from \cite{GurjarMiyanishi:affinesurfaceslkd1}
(see also \cite[Theorem 2.6]{Zaidenberg:1998exot}, \cite{tomDieckPetrie:contractibleaffinesurfaces} and
\cite{FlennerZaidenberg:contractiblelkd1}).
We start with a line arrangement in $\P^2$ as in Figure \ref{fig:linearrangement}.

\begin{tikzpicture} [scale=1.0] \label{fig:linearrangement}
%

\draw (5,5) -- (5,0);
\draw (5,8) -- (5,5);
\draw (5,5) circle (0.05);
\draw (3,8) -- (3+2*8/3,0);
\draw (7,8) -- (7-2*8/3,0);
\draw (0,2) -- (10,2);
\draw[fill=black]  (7-2*4.5/3,3.5) node[right] {$p_0$} circle (0.05);
\draw[fill=black]  (5,2) circle (0.05);
\draw[fill=black]  (3+2*6/3,2) circle (0.05);
\node at (6,1.8) {$\cdots$};
\node[right] at (5,1.8) {$p_1$};
\node[left] at (3+2*6.2/3,1.8) {$p_s$};
\node[left] at (7-2*4/3,4) {$D_0$};
\node[right] at (5,4) {$D_1$};
\node[right] at (3+2*4/3,4) {$D_s$};
\node[right] at (10,2) {$H$};

\end{tikzpicture}

Here we have curves $D_0,D_1,\cdots,D_s,H$ in this line arrangement.
Let $W$ be the divisor in $\P^2$ representing this line arrangement.
At the point $p_0$ we blow up our surface many times according to the following rules:
\begin{enumerate}
\item The first blow up must be at $p_0$.
\item Each subsequent blow up must be
on the exceptional divisor of the previous blow up and at a smooth point
of the total transform of $W$.
\end{enumerate}
At the point $p_i$ where $i > 0$, we blow up in such a way as to to resolve
the point of indeterminacy of $\frac{x^{m_i}}{y^{n_i}}$
(viewed as a birational map to $\P^1$) where $x,y$
are local coordinates around $p_i$ with $W = \{xy = 0\}$.
These are chains of blowups where we only blow up along points $p$ satisfying:
\begin{enumerate}
\item $p$ is in the exceptional divisor of the previous blowup.
\item $p$ is a nodal singular point of the total transform of $W$.
\end{enumerate}
Let $E_0,E_1,\cdots,E_s$ be the last exceptional curves in these chains of blowups
over $p_0,\cdots,p_s$.
We let $X$ be equal to $\P^2$ blown up as described above
and we let $D$ be the divisor in $X$ equal to the total transform of $W$
minus the last exceptional curves $E_0,\cdots,E_s$.
Our surface $A$ is equal to $X \setminus D$.
The integers $m_i, n_i$ are coprime and satisfy a certain equation
to ensure that our surface $A$ is affine and acyclic of log Kodaira dimension $1$.

\begin{lemma} \label{lemms:unirulednessforlkd1}
Suppose $A$ is an acyclic surface of log Kodaira dimension $1$.
Then it is compactified $2$ uniruled.
\end{lemma}
\proof of Lemma \ref{lemms:unirulednessforlkd1}.
We will use the notation $X,D_0,\cdots,D_s,E_0,\cdots,E_s$, $H$, $p_0,\cdots,p_s$,
from before this lemma to describe $A$.
 We have three cases:
\begin{enumerate}
\item \label{item:caseblowuponce}
$s = 1$.
\item \label{item:casemorefibers1}
$s > 1$ and $E_i$ intersects $H$ for some $i$.
\item \label{item:casemorefibers2}
$s > 1$ and $E_i$ does not intersect $H$ for any $i$.
\end{enumerate}

{\it Case (\ref{item:caseblowuponce}):}
Let $\text{Bl}_{p_0}(\P^2)$ be the blowup of $\P^2$ at the point $p_0$.
We have a map from $X$ to $\text{Bl}_{p_0}(\P^2)$ which is a sequence of blowdown maps.
We also have a fibration $\pi : \text{Bl}_{p_0} \to \P^1$
whose fibers are proper transforms of lines in $\P^2$ passing through $p_0$.
Let $\widetilde{\pi} : X \to \P^1$ be the composition of the map from $X$ to $\text{Bl}_{p_0}(\P^2)$ with
$\pi$.
Because $s = 1$, we have that a generic fiber of $\widetilde{\pi}$
intersects $D$ twice.

Also, the proper transform of $D_0$ is a fiber of $\widetilde{\pi}$.
So if we choose any almost complex structure $J$ which is equal to the standard
one near $D$, we have a fiber represented by the irreducible $J$ holomorphic curve $D_0$.
So, $D_0$ has intersection number $0$ with any fiber.
By Corollary \ref{corollary:uniruledfibrationproperty}
we then get that $A$ is compactified $2$ uniruled.

{\it Case (\ref{item:casemorefibers1}):}

Let $q$ be the point where all the divisors $D_0,\cdots,D_s$ intersect in one point
and let $\widetilde{q}$ be the corresponding point in $X$.
We blow up $\widetilde{q}$ giving us $\widetilde{X}$.
Let $\widetilde{D}$ be the total transform of $D$, so $A = \widetilde{X} \setminus \widetilde{D}$.
Note that $\widetilde{X}$ is equal to $\text{Bl}_q \P^2$ blown up many times at the points
$p_0,\cdots,p_s$. Hence we have a natural blowdown map
$\text{Bl}_{\widetilde{X}}: \widetilde{X} \to \text{Bl}_q \P^2$.
There is a natural map $\pi' : \text{Bl}_q \P_2 \to \P^1$ whose fibers are proper transforms of lines passing through $q$.
Let $\widetilde{\pi}' : \widetilde{X} \to \P^1$ be the composition
$\pi' \circ \text{Bl}_{\widetilde{X}}$.
We let $\widetilde{E}_0,\cdots,\widetilde{E}_s$ be the proper transforms of
$E_0,\cdots,E_s$ in $\widetilde{X}$ respectively.
We similarly define $\widetilde{H}$, $\widetilde{D}_i$.
Let $E$ be the proper transform of the exceptional divisor of $\text{Bl}_q \P^2$
in $\widetilde{X}$.
The image of our morphism $\widetilde{\pi}'$ is $B := \P^1$.
We define $b_j \in B$ so that $(\widetilde{\pi}')^{-1}(b_j)$ contains $\widetilde{E}_j$.

We have that $\widetilde{E}_i$ intersects $\widetilde{H}$ for some $i$. This is contained in some fiber
$(\widetilde{\pi}')^{-1}(b_i)$.
We have that $(\widetilde{\pi}')^{-1}(b_i)$ is obtained from $(\pi')^{-1}(b_i)$
by blowing up the point where this fiber intersects $H$ repeatedly.
Hence if $R$ is the irreducible component of $(\widetilde{\pi}')^{-1}(b_i)$ that intersects $E$, then it is smooth and has self intersection $-1$.
This means that $R + E$ is an effective nef divisor.
We have that $\widetilde{E}_0.D = 1$
so let $D_E$ be the unique divisor that intersects $\widetilde{E}_0$.
Let $D'$ be the union of irreducible curves in $D$ not intersecting $\widetilde{E}_0$
and let $\Delta'$ be the connected component of $D'$ containing $R \cup E$.
Using Lemma \ref{lemma:nefdivisor} with the divisors $\Delta' + D_E$ and $E$ there is a nef divisor $G$
with the property that $G.D_E \neq 0$ and $G.E = 0$.
The generic fiber  of $\widetilde{\pi}'$ intersects $\widetilde{D}$ twice.
Also $D_E \cup \widetilde{E}_0$ is contained in $(\widetilde{\pi}')^{-1}(b_0)$
and also there is an effective divisor $D_F$ whose support is in $(\widetilde{\pi}')^{-1}(b_0) \cap D$
and $\kappa \in \N$ so that $[D_F] + \kappa [\widetilde{E}_0]$ is homologous to a fiber of $\widetilde{\pi}'$.
So by Lemma \ref{lemma:dimension2uniruled4} we get that $A$ is compactified $2$ uniruled.

{\it Case (\ref{item:casemorefibers2}):}

Because $s>1$ and $E_i$ does not intersect $H$ for all $i$, we get that $E_1$ and $E_2$ exist and do not intersect $H$.
We have that $(\widetilde{\pi}')^{-1}(b_j)$ is a union of irreducible curves
$F_j$ in $D$ plus $\widetilde{E}_j$.
Also $\widetilde{E}_j . \widetilde{D} = \widetilde{E}_j.\widetilde{F}_j$ for all $j$
and $\widetilde{D} \cup \widetilde{E}_1 \cup \widetilde{E}_2$ is a smooth normal crossing divisor.
Let $D'_j$ be equal to $\widetilde{D}$ minus the irreducible components of $\widetilde{D}$ intersecting $E_j$.
We have that $D'_j$ is connected for each $j$ and every irreducible component of $\widetilde{D}$ which intersects $E_j$
also intersects $D'_j$.

We have that $\widetilde{D}$ is connected and
$E \cup \widetilde{D}_0$ are disjoint from $E_j$ for each $j$.
Both $E$ and $\widetilde{D}_0$ intersect each other and have self intersection
$-1$  which implies that $G := E + \widetilde{D}_0$ is nef and contained in $D'_j$ for each $j$.
This does not intersect $E_j$.
Hence by Lemma \ref{lemma:dimension2uniruled3} we then get that
$A$ is compactified $2$ uniruled.
\qed


\proof of Theorem \ref{theorem:acyclicsurfaceinvariance}.
In order to prove this theorem we only need to show the following fact:
{\it If $A$ has log Kodaira dimension $i$ where $i \leq 1$, then
$B$ has log Kodaira dimension $ \leq i$.}
This is because log Kodaira dimension is at most $2$.

By \cite{Fujita:1982alg} we have that the log Kodaira dimension of $A$
is either $-\infty$, $1$ or $2$. Also if it is equal to $-\infty$ then
$A = \C^2$.
Suppose the log Kodaira dimension of $A$ is $-\infty$ then
$A = \C^2$. Also $B$ is diffeomorphic to $A$ and hence contractible and simply connected at infinity.
By \cite{Ramanujam:affineplane} we then get that $B$ is isomorphic to $\C^2$
and hence has log Kodaira dimension $-\infty$.
Now suppose that $A$ has log Kodaira dimension $1$.
By Lemma \ref{lemms:unirulednessforlkd1} we have that $A$ is compactified $2$ uniruled,
so by Lemma \ref{lemma:compactifieduniruledlkd}, $B$ has log Kodaira dimension $\leq 1$.
Putting everything together gives us that $A$ and $B$ must have the same log Kodaira dimension.
\qed

\subsection{Dimension 3}

\begin{theorem} \label{theorem:logkodairaresultsindimension3}
Suppose that $A$ is a smooth affine variety of dimension $3$ such that $A$
admits a compactification $X$ with the following properties:
\begin{enumerate}
\item The compactification divisor $D$ is smooth normal crossing and nef.
\item The linear system $|D|$ contains a smooth member.
\end{enumerate}
Let $B$ be any smooth affine variety symplectomorphic to $A$ and
$\overline{\kappa}(A) = 2$ then $\overline{\kappa}(B) \leq 2$.
\end{theorem}
\proof of Theorem \ref{theorem:logkodairaresultsindimension3}.
By \cite{Kishimoto:affinethreefolds} we have that $A$ admits a $\C^*$ fibration.
In fact we can say more:
\cite[Lemma 4.1,4.2,4.3]{Kishimoto:affinethreefolds} says that there is a projective variety
$X^{s}$ and a nef divisor $D^{s}$ so that
\begin{enumerate}
\item  $A = X^{s} \setminus D^{s}$.
\item There is a morphism $\pi : X^{s} \to W$ of projective varieties
with the property that a generic fiber is isomorphic to $\P^1$ and intersects $D^{s}$ twice.
\end{enumerate}
By \cite{hironaka:resolution} we can blow up $X^{s}$ away from $A$ giving us a smooth projective
variety $X$ and so that the total transform $D$ of $D^{s}$ is a smooth normal crossing divisor.
Let $\widetilde{\pi} : X \to W$ be the composition of $\pi$ with the blowdown map
$X \to X^{s}$.
Let $D_X$ be the effective divisor which is the pullback of $D^{s}$ under the blowdown map.
We have that $D_X$ is nef  and that a generic fiber $F$ of $\widetilde{\pi}$
satisfies $F.D_X = 2$.
By Lemma \ref{lemma:nefdivisoruniruled}, we then get that $A$ is compactified $2$ uniruled.
Hence by Lemma \ref{lemma:compactifieduniruledlkd}, we get that
the log Kodaira dimension of $B$ is $\leq 2$.
\qed

\section{Uniruledness of compactifications} \label{section:unirulednesscompactifications}

If a projective variety $X$ has a morphism $f : \P^1 \rightarrow X$ passing through every point $x \in X$
then we say that $X$ is {\it uniruled}.

\begin{theorem} \label{theorem:birationalinvarianceofuniruledness}
Suppose that two smooth projective varieties $P$ and $Q$ have affine open subsets
$A$, $B$ with the property that $A$ is symplectomorphic to $B$.
Then $P$ is uniruled if and only if $Q$ is.
\end{theorem}
\proof of Theorem \ref{theorem:birationalinvarianceofuniruledness}.
Suppose that $P$ is uniruled, then we will show that $Q$ is uniruled.
Let $D_P$ be an effective ample divisor whose support is $P \setminus A$
and $D_Q$ an effective ample divisor whose support is $Q \setminus B$.
By \cite{Ruan:virtual} or \cite{Kollar:lowdegree} we have that
$\langle [\text{pt}], \alpha_1, \cdots, \alpha_k \rangle^P_{0,\beta} \neq 0$
for some $\beta \in H_2(P,\Z)$ and cohomology classes $\alpha_1,\cdots,\alpha_k$.
Let $k := \beta. D_P$.
This means that any compatible $J$ in $P$ which is standard near $D_P$ has the property that
there is some $J$ holomorphic curve $u : \Sigma \to P$ passing through any point $p$.
Because $D_P$ is nef, each irreducible component $\Sigma_i$ of $\Sigma$
satisfies $u_*(\Sigma_i).D_P \leq k$. In particular this is true for
any irreducible component that passes through $p$.
Hence $A$ is compactified $k$-uniruled.
So by Theorem \ref{theorem:algebraicuniruledimplesuniruled},
we have that the Liouville domain $\overline{A}$ is $(k,\Lambda)$-uniruled
for some $\Lambda > 0$.
Because the completion of $\overline{A}$ is symplectomorphic to the completion of $\overline{B}$
we then get by Theorem \ref{theorem:symplecticinvariance} that
$\overline{B}$ is $(k,\Lambda')$ uniruled for some $\Lambda' > 0$.
Hence by Theorem \ref{theorem:kuniruledimpliesalgebraicallyuniruled},
we have that $B$ is algebraically $k$-uniruled.
This implies that its compactification $Q$ is uniruled.
By symmetry, if $Q$ is uniruled then $P$ is. Hence $P$ is uniruled if and only if $Q$ is uniruled.
\qed

\section{Appendix : plurisubharmonic functions on smooth affine varieties}

The contents of this appendix are all contained inside the proof
of \cite[Lemma 2.1]{McLean:affinegrowth} and the ideas of that proof are contained
in \cite[Section 4b]{Seidel:biasedview}.
We let $A$ be a smooth affine variety.
Here we recall the construction of the Liouville domain $\overline{A}$
(see Definition \ref{defn:canonicasymplecticformassociatedtoA}).
Choose any algebraic embedding $\iota$ of $A$ into $\C^N$
(so it is a closed subvariety).
We have $\theta_A := -d^c R = \sum_i\frac{r_i^2}{2} d\vartheta_i$
where $(r_i,\vartheta_i)$ are polar coordinates for the $i$th $\C$ factor.
We have that $d\theta_A$ is equal to the standard symplectic structure on $\C^N$.
By abuse of notation we write $\theta_A$ for $\iota^* \theta_A$,
and $\omega_A := d\theta_A$. Here $(\overline{A},\theta_A)
:= (R^{-1}(-\infty,C],\theta_A)$ for $C \gg 0$.

We can also construct other Liouville domains as follows
(see Definition \ref{defn:symplecticformonaffinevariety}):
Let $X$ be a smooth projective variety
such that $X \setminus A$ is a smooth normal crossing divisor
(an SNC compactification).
Let $L$ be an ample line bundle on $X$ given by an effective divisor $D$ whose support is $X \setminus A$.
From now on such a line bundle will be called a {\bf line bundle associated to an SNC compactification $X$ of $A$}.
Suppose $|\cdot|$ is some metric on $L$
whose curvature form is a positive $(1,1)$ form. Then if $s$ is some section of $L$ such that $s^{-1}(0) = D$
then we define $\phi_{s,|\cdot|} := -\log{|s|}$ and $\theta_{s,|\cdot|} := -d^c \phi_{s,|\cdot|}$. The two form $d\theta_{s,|\cdot|}$
extends to a symplectic form $\omega_{|\cdot|}$ on $X$ (which is independent of $s$ but does depend on $|\cdot|$).
We will say that $\phi_{s,|\cdot|}$ is a {\it plurisubharmonic function associated to $L$},
$\theta_{s,|\cdot|}$ a {\bf Liouville form associated $L$} and
$\omega_{|\cdot|}$ a {\bf symplectic form on $X$ associated to $L$}.
From \cite[Section 4b]{Seidel:biasedview}, we have that for
$C \gg 1$,\[(A_C,\theta_C) := (\phi_{s,|\cdot|}^{-1}(-\infty,C], \theta_{s,|\cdot|})\]
is a Liouville domain.

Let $(r_i,\vartheta_i)$ be polar coordinates for the $i$'th
factor in $\C^N$.
\begin{lemma} \label{lemma:nosingularitiesofaparticularfunction}
If we compactify $\C^N$ by $\P^N$, there is a section $S$ of ${\mathcal O}(1)$
and metric $\|\cdot\|$ with the following properties:
\begin{enumerate}
\item $-\log{\|S\|}|_A$ is equal to $f(R)$ for some non-decreasing smooth function $f : \R \to \R$.
\item  $-\log{\|S\|}|_A$ has no singularities near infinity.
\end{enumerate}
Hence $R|_B$ has no singularities near infinity.
\end{lemma}
\proof of Lemma \ref{lemma:nosingularitiesofaparticularfunction}.
Let $H := \P^N \setminus \C^N$ and
let $S$ be a section of ${\mathcal O}(1)$ such that $S^{-1}(0) = H$.
Let $\|\cdot\|$ be the standard Fubini Study metric on ${\mathcal O}(1)$.
We have that $U(N+1)$ acts on $\P^N$
and it naturally lifts to an action on the total space of ${\mathcal O}(1)$.
Let $U(N) \subset U(N+1)$ be the natural subgroup that preserves $H$.
We have that $\|S\|$ is invariant under this action.

Because $-\log{\|S\|}$ is invariant under this action and exhausting, it is equal to
$f(R)$ for some non decreasing smooth function $f : \R \to \R$.
This is because $U(N)$ acts transitively on the level sets of $R$.

Let $X$ be the closure of $A$ in $\P^N$.
By \cite{hironaka:resolution} we can blow up $\P^N$ along $H$
so that the proper transform $\widetilde{X}$ of $X$ is smooth
and the total transform $\widehat{H}$
of $H \cap X$ inside $\widetilde{X}$ is a smooth normal crossing divisor.
We pull back the line bundle ${\mathcal O}(1)|_X$ to a line bundle $L_X$ on $\widetilde{X}$
and also pull back the metric $\|\cdot\|$ and section $S$.
We will write $\|\cdot\|_X$ and $S_X$ for the new metric and section.

Let $p \in \widehat{H}$ and choose local holomorphic coordinates
$z_1,\cdots,z_n$ on $\widetilde{X}$ and a trivialization of $L_X$ around $p$ so that
$S_X = z_1^{w_1} \cdots z_n^{w_n}$ ($w_i \geq 0$).
The metric $\|.\|_X$ on $L_X$ is equal to $e^{\psi}|.|$
for some function $\psi$ with respect to this trivialization
where $|.|$ is the standard metric on $\C$.
So \[-d\log{\|S_X\|_X} = -d\psi - (\sum_i w_i d\log{|z_i|}).\]
If we take the vector field
$Y := -r_1 \partial_{r_1} \cdots - r_n \partial_{r_n}$ (where $z_j = r_je^{i\vartheta_j}$),
then $d\log{(|z_j|)}(Y) = -1$ and $d\psi(Y)$ tends to zero.
Hence $d\log{\|S_X\|_X}$ is non-zero near infinity which implies that
$f(R)|_A = -\log{\|S\|_X}$ has no singularities near infinity.
\qed

\begin{lemma} \label{lemma:affineLiouvilledomain}
For $C \gg 1$, $(A_C,\theta_C)$ is Liouville deformation equivalent to
$(\overline{A},\theta_A)$.
\end{lemma}
\proof of Lemma \ref{lemma:affineLiouvilledomain}.
By Lemma \ref{lemma:nosingularitiesofaparticularfunction}
we have that $R|_A$ has no singularities for $R \geq C$.
Let $c \geq C$ and write  $A'_c := (R|_A)^{-1}(-\infty,c]$.
Because $c$ is a regular value of $R|_A$,
we have that $A'_c$ is a Liouville domain and by definition it is equal to
$(\overline{A},\theta_A)$.

Let $S$ and $\|\cdot\|$ be the section and metric on ${\mathcal O}(1)$
coming from Lemma \ref{lemma:nosingularitiesofaparticularfunction}.
We also have that $f(R) = -\log{\|S\|}$ where $f$ is a smooth function with
positive derivative when $R$ is large.
Let $A''_c := (-\log{\|S\|})^{-1}(-\infty,c]$.
We have that $A'_{f(c)} = A''_c$ for $c \gg 1$.
We also have that $t \theta_A + (1-t) \theta_{S,\|\cdot\|}|_A$
is a deformation of Liouville domains from
$(A'_{f(c)},\theta_A)$ to $(A''_c, \theta_{S,\|\cdot\|})$.

Let $\phi_{s,|\cdot|}$ be a plurisubharmonic function associated to our line bundle $L$
so that $A_C = \phi_{s,|\cdot|}^{-1}(-\infty,C]$.
Let $\phi_t := (1-t)\phi_{s,|\cdot|} - t \log{\|S\|}$
and let $A^t_c := \phi_t^{-1}(-\infty,c]$.
By using work from \cite[Section 4b]{Seidel:biasedview},
we have for $C$ large enough that
$(A^t_C, -d^c \phi_t)$ is a Liouville deformation
from $(A''_C,\theta_{S,\|\cdot\|})$ to
$(A_C,\theta_C)$.
Hence by composing the above two Liouville deformations, we get that
$(A_C,\theta_C)$ is Liouville deformation equivalent to
$(\overline{A},\theta_A)$ for $C$ large enough.
\qed

\bibliography{references}

\end{document}